\documentclass[11pt]{amsart}
\usepackage{amsmath,amsfonts,amssymb,amsxtra,latexsym,amscd,enumerate,amsthm}

\usepackage{latexsym}
\usepackage{epsfig}
\usepackage{verbatim}

\def\r{\mathcal{R}}
\def\g{\gamma}
\def\a{\alpha}

\def\d{\delta}

\def\b{\beta}

\def\p{\Phi}
\def\f{\varphi}
\def\ep{\epsilon}
\def\es{\emptyset}
\def\l{\lambda}
\def\o{\omega}

\def\R{\mathbb{R}}
\def\C{\mathcal{C}}
\def\Z{\mathbb{Z}}
\def\K{\mathcal{K}}
\def\D{\mathcal{D}}
\def\P{\mathbb{P}}
\def\p{\mathcal{P}}
\def\l{\mathcal{L}}

\def\J{\mathcal{J}}
\def\M{\mathcal{M}}
\def\T{\mathcal{T}}
\def\A{\mathcal{A}}
\def\B{\mathcal{B}}
\def\TT{\mathbb{T}}
\def\N{\mathbb{N}}
\def\Z{\mathbb{Z}}
\def\beq{\begin{equation}}
\def\eeq{\end{equation}}

\def\beq{\begin{equation}}
\def\eeq{\end{equation}}

\newtheorem{d0}{Definition}
\newtheorem{c0}{Corollary}
\newtheorem{t1}{Theorem}
\newtheorem{l1}{Lemma}
\newtheorem{p1}{Proposition}
\newtheorem*{claim}{CLAIM}

\begin{document}
\title[]{The (weak-$L^2$) Boundedness of The Quadratic Carleson Operator}

\author{Victor Lie}

\date{\today}

\address{Department of Mathematics, UCLA, Los Angeles CA 90095-1555}

\email{vilie@math.ucla.edu}
\address{Institute of Mathematics of the Romanian Academy, Bucharest, RO
70700 \newline \indent  P.O. Box 1-764}


\maketitle

\begin{abstract}

We prove that the generalized Carleson operator with polynomial
phase function of degree two is of weak type (2,2). For this, we
introduce a new approach to the time-frequency analysis of the
quadratic phase.
\end{abstract}
$\newline$
\section{\bf Introduction}

 The historical motivation for the subject of this paper
is rooted in Luzin's Conjecture (1913), which asserts that the
Fourier series of a function $f\in L^2(\TT)$ converges pointwise to
$f$ Lebesgue-almost everywhere. In 1966, L. Carleson gave a positive
answer to this conjecture in the celebrated paper \cite{1}. His
result can be essentially reformulated - via \cite{8} - as follows:
$\newline$ {\bf Theorem 0.} {\it If for $f\in C^{1}(\TT)$ we define
the expression\footnote{In what follows, we will always omit
principal value notation.} \beq\label{car}
Cf(x):=\sup_{a>0}\left|\int_{\TT}\frac{1}{y}e^{iay}f(x-y)dy\right|,
\eeq then $C$ is of weak type (2,2), i.e.: \beq\label{v2}
\left\|Cf\right\|_{L^{2,\infty}(\TT)}\leq A
\left\|f\right\|_{L^2(\TT)}\:, \eeq where here, by convention,
$\TT=[-\frac{1}{2},\frac{1}{2}]$ and $A>0$ is an absolute constant.}

In addition to Carleson's proof, we point out two more proofs of the
above result: one due to Fefferman \cite{2}, using a very beautiful
geometric combinatorial argument, and the other due to Lacey and
Thiele \cite{5}, inspired from the subtle techniques they developed
for proving the Calder\'on conjecture (\cite{6} and \cite{7}). Now,
given the statement of Theorem 0, it is natural to hope that this
result may be set in a broader context. Following this direction,
Stein conjectured that the generalized Carleson operator defined by
\beq\label{v3}
C_{d}f(x):=\sup_{deg(P)=d}\left|\int_{\TT}\frac{1}{y}e^{iP(y)}f(x-y)dy\right|
\eeq (here $d\in \Z\:\:,\:\:d\geq 2$, $P$ is a polynomial of degree
$d$, and $f\in C^{1}(\TT)$) $\newline$ obeys the same bounds as $C$.

In \cite{9} he proved this conjecture, subject to the key
restriction that the supremum in \eqref{v3} be taken in the class of
quadratic polynomials with no linear term. Further, using the $TT^*$
method and a variant of van der Corput's lemma, Stein and Wainger
\cite{10} extended this result for polynomials of any degree, but again
without the first degree term.

Our aim in this paper is to provide a positive answer to this
conjecture for the case $d=2$.\footnote{We mention here that (using
similar techniques to those in \cite{5}) M. Lacey published in
\cite{4} a proof of this result, but as was revealed later by A.
Ionescu, this was incorrect - for details see \cite{3}.}

The main result of the article is given by:

\begin{t1}\label{ma}
Let $1\leq p < 2$; then the expression \beq\label{v4} Tf(x):=
C_{2}f(x)=\sup_{a,b\in
\R}\left|\int_{\TT}e^{i\left\{ay\:+\:by^2\right\}}\frac{1}{y}f(x-y)dy\right|
\eeq satisfies
$$\left\|Tf\right\|_{L^{p}(\TT)}\lesssim_{p}\left\|f\right\|_{L^{2}(\TT)}.$$
\end{t1}

Combining this result with the techniques developed by Stein in
\cite{8}, we easily deduce:
\begin{c0}\label{co}
$T$ is of weak type (2,2).
\end{c0}

The proof of Theorem 1 is a combination of analytic and geometric
facts; it relies on a new perspective of the time-frequency
localization of the quadratic phase to which we adapt the techniques
presented in \cite{2}.

One particular feature of this paper is that it presents for the
first time a time-frequency proof of the boundedness of a maximal
operator which is invariant under quadratic modulations.

Another novelty of this paper is that we show that one can prove the
(Quadratic) Carleson Theorem using a single dyadic grid partition
(on each axis defining the time-frequency plane).\footnote{For this,
we will involve in our reasoning certain dilation factors of our
tiles (see Section 7 - forest decomposition algorithm).}

Finally, given the powerful geometric intuition developed in
Fefferman's paper, and also the fact that many of the reasonings
here rely on his work, we have chosen to present our paper
maintaining the structure and some of the notations appearing in
\cite{2}.

{\bf Acknowledgements.} I would like to thank Ciprian Demeter, Camil
Muscalu and Terence Tao for reading the paper and offering several
useful comments, and Zubin Gautam for his care in improving the
English presentation. I am grateful to my advisor Christoph Thiele
for suggesting this problem and for useful discussions. Finally, I
am indebted to Nicolae Popa  for introducing me to the field of
harmonic analysis and for offering his constant moral support.

\section{\bf Preliminaries and outline of the proof}

As our problem is of a time-frequency nature, it will be based on
two steps:

(A) - a discretization procedure, in which we split our operator
into ``small pieces" that are well-localized in both time and
frequency.

(B) - a selection algorithm, which relies on finding (qualitative
and quantitative) criteria depending on which we decide how to glue
the above-mentioned pieces together to obtain a global estimate on
our operator.

For task (A), we first need to study the symmetries of our operator.
This is because these symmetries will determine the geometric
properties of the time-frequency portrait of our operator,
properties that will provide a significant indication of how to
naturally decompose the operator ``along its fibers."

We define the following classes of symmetries\footnote{Since the symmetries 1) and 2) do not preserve the periodicity of the object on which they are acting, in what follows one should regard $L^2(\TT)$ as the space of functions which are $L^2$-integrable  on any given unit interval.}:
$\newline$1) Modulations:
 $$M_{a}:\:\:\:L^2(\TT)\longrightarrow L^2(\TT)\:\:(a\in \R)\:\:\:\:\textrm{by}\:\:\:\:M_{a}f(x):=e^{iax}\:f(x)$$
$\newline$2) Quadratic Modulations:
 $$Q_{b}:\:\:\:L^2(\TT)\longrightarrow L^2(\TT)\:\:(b\in \R)\:\:\:\:\textrm{by}\:\:\:\:Q_{b}f(x):=e^{ibx^2}\:f(x)$$
$\newline$3) Translations:
 $$\tau_{y}:\:\:\:L^2(\TT)\longrightarrow L^2(\TT)\:\:(y\in \R)\:\:\:\:\textrm{by}\:\:\:\:\tau_{y}f(x):=f(x-y)$$
$\newline$4) Dilations:
 $$D_{\lambda}:\:\:\:L^2(\TT)\longrightarrow L^2(\TT)\:\:(\lambda\in \N)\:\:\:\:\textrm{by}\:\:\:\:D_{\lambda}f(x):=f(\lambda x)$$

The key observation is that we can recover the operators $C$ and $T$
from the action of these symmetries (particularly 1) and 2)) on the
Hilbert transform\footnote{Strictly speaking, the kernel of the
Hilbert transform should be $\cot{\pi y}$; for convenience, we work
instead with $\frac{1}{y}$.}, defined by
 $$H:\:\:\:L^2(\TT)\longrightarrow L^2(\TT)\:\:\:\:\:\:\:\:\:\:Hf(x):=\int_{\TT}\frac{1}{y}f(x-y)dy\:.$$

Indeed, the periodic Carleson operator

$$Cf(x)=\sup_{a\in\: \R}\left|\int_{\TT}\frac{1}{y}\:e^{iay}f(x-y)dy\right|$$
can be rewritten as \beq\label{v5}
 Cf(x)=\sup_{c\in\: \R}\left|M_{c}\:H\:M_{c}^{*}f(x)\right|,
 \eeq
 while our periodic Quadratic Carleson operator
 $$Tf(x)=\sup_{a,b\in\: \R}\left|\int_{\TT}e^{i\left\{ay\:+\:by^2\right\}}\frac{1}{y}f(x-y)dy\right|$$
can be rewritten as \beq\label{v6} Tf(x)=\sup_{b,c\in\:
\R}\left|M_{c}Q_{b}\:H\:Q_{b}^{*}M_{c}^{*}f(x)\right|. \eeq (Remark
that in the previous formulas the action of translations and
dilations is hidden in the structure of the Hilbert transform, which
is the unique - up to identity - $L^2$-bounded linear operator that
commutes with both symmetries.)

 These facts help us to conclude that $C$ essentially\footnote{The relations 1)-4) are literally true if we work in the setting of $\R$ rather than $\TT$.  Relations 1), 3), and 4) remain true in the torus case, while 2) serves as a useful heuristic (especially for the operator $T$) inherited from the real case.} obeys the relations
 $\newline$ 1) $C \tau_y = \tau_y C$
 $\newline$ 2) $C D_{\lambda} = D_{\lambda}C$
 $\newline$ 3) $C M_a = C$
 $\newline$
while for the operator $T$, besides the analogous relations we have the
extra condition $\newline$ 4) $T Q_{b} =T\: .$

We now analyze the effect of these symmetries on the time-frequency
decomposition of our operator $T$.  To help build up intuition, we
will consider three cases of increasing complexity: the Hilbert
transform $H$, the Carleson operator $C$, and finally the Quadratic
Carleson operator $T$.

As announced, we first look at the simplest object, namely the
Hilbert transform; we begin by isolating the kernel and splitting
 it - taking advantage of the dilation symmetry of $H$ - as follows:
 $$\frac{1}{y} = \sum_{k\in \N}\psi_k (y)$$
where $\psi\in C_{0}^{\infty}$ is an odd function supported away
from the origin and  $\psi_k(y)=2^{k}\psi(2^{k}y)$, $k\in \N$;
consequently,
 \beq\label{hilb}
 Hf(x)\: = \:\sum _{k\in \N} \int\psi_k(y)f(x-y)dy\:.
 \eeq

Now for each scale $k$ we take the collection $\{I_{k,j}\}_j$ of all
dyadic intervals in $[0,1]$ of length $2^{-k}$. Using the
translation invariant property of $H$ we write \beq\label{hils}
  Hf(x)\: =
  \:\sum_{k,j}H_{k,j}f(x)\:=\:\sum_{k,j}\left\{\int\psi_k(y)f(x-y)dy\right\}\chi_{I_{k,j}}(x)\:,
  \eeq where $\chi_I$ is, as usual, the characteristic function of $I$.

Now each $H_{k,j}f$ has time support included in $I_{k,j}$ while on
the frequency side it is ``morally" supported near the origin, in an
interval of length $|I_{k,j}|^{-1}$. Consequently, the
time-frequency picture of $H$ is as given in Figure 1.

The above story can  be expressed more intuitively as follows:
Observe that the translation symmetry acts on the $j$-direction,
while the dilation symmetry acts on the $k$-direction.  If we
approximate the piece $H_{1,1}f$ by a smooth compactly supported
function $\varphi_0$, then the time-frequency portrait of
$\varphi_0$ is a square of area one located near the origin.  Since
$Hf$ is, roughly speaking, just a sum of dilations and translations
of $\varphi_0$, by basic properties of the Fourier transform we
obtain Figure 1 as the time-frequency picture of $Hf$. From the
figure we also note that the origin plays a special role in this
decomposition.
\begin{figure}[!h]
\begin{center}
\epsfig{file=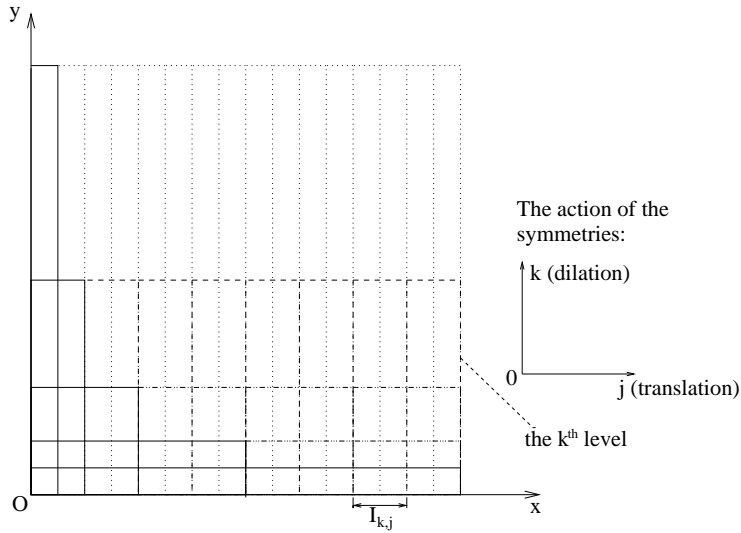,height=7cm}
\end{center} \caption{The time-frequency decomposition of the Hilbert transform}
\end{figure}\\

We now consider the Carleson operator  as described in \eqref{v5}.
In this case we will have to deal with one more symmetry given by
the modulation invariance property, so we will try first to
understand a simpler situation, namely how $M_c $ acts on a smooth
compactly supported function $\f$. As we may remark from Figure 2,
in the time-frequency plane, $M_c$ will translate the rectangle
representing the localization of $\f$ by $c$ units in the frequency
direction.
\begin{figure}[!h]
\begin{center}
\epsfig{file=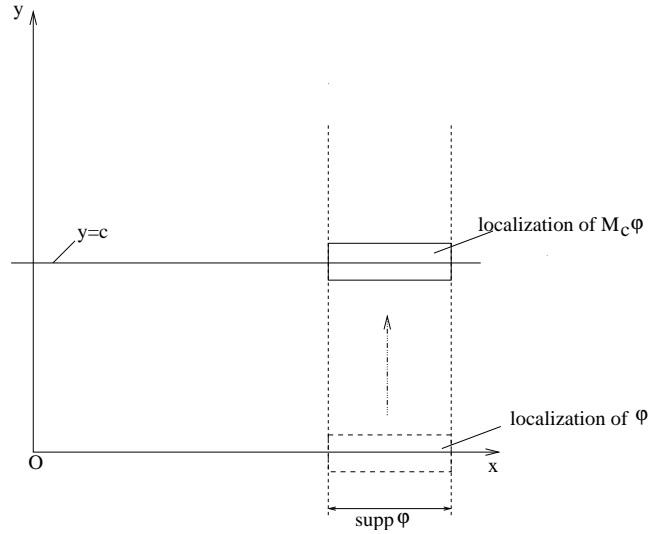,height=7cm}
\end{center} \caption{The time-frequency localization of $M_c\f$}
\end{figure}\\

Now, as in  \eqref{hilb}, we have that
$$M_{c}HM_{c}^{*}f(x)\:=\:\sum_{k,j}M_{c}H_{k,j}M_{c}^{*}f(x)\:=\:\:\sum_{k,j}\left\{\int(M_{c}\psi_k)(y)f(x-y)dy\right\}\chi_{I_{k,j}}(x)\:,$$
and so combining this with the previous observation, we deduce that
the time-frequency picture of $M_c H M_c^{*} $ will be nothing more
than a frequency-translation by $c$ units of the corresponding
picture of $H$.

Exploiting this fact in the form of \eqref{v5}, we conclude that the
time-frequency localization of $C$ is as presented in Figure 3.
\begin{figure}[!h]
\begin{center}
\epsfig{file=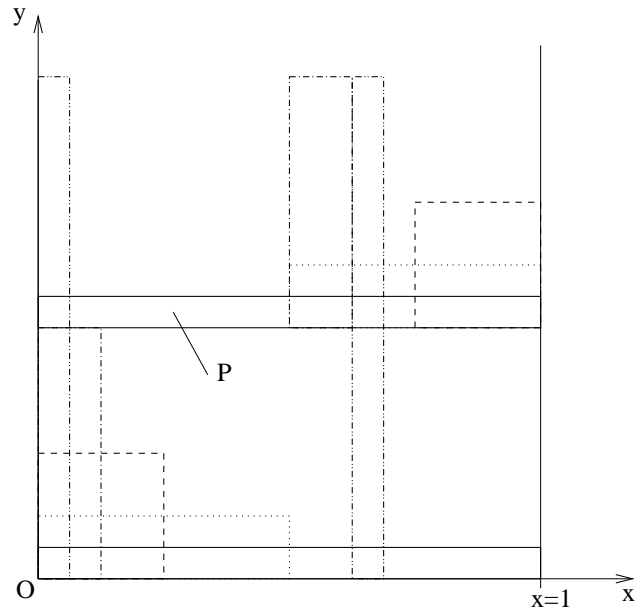,height=8cm}
\end{center} \caption{The time-frequency decomposition of the Carleson Operator}
\end{figure}\\

Remark that, unlike the Hilbert transform case, there is no
preferential point in the splitting of $C$. Also, this picture
suggests that $C$ may be written (after a linearization procedure)
as
  $$Cf=\sum_{P}C_Pf$$
with each $C_P$ a linear operator localized in a certain
(Heisenberg) rectangle $P$ (Figure 3). This is a key observation
used explicitly in both \cite{2} and \cite{5}.

We finalize part (A) of our program with the analysis of our
operator $T$. As before, we begin by isolating the extra symmetry -
$ Q_b$, that adds to those appearing in the previous cases. We will
approach the study of the time-frequency representation of this
quadratic symmetry from two perspectives: ({\it q1}) a {\it
restrictive} one and ({\it q2}) a {\it relational} one.

({\it q1}) The {\it restrictive} perspective relies on the following
basic approach: given an object (Schwartz function on $\R$) - call
it $h$ - describe (in terms of a picture) the space and frequency
regions\footnote{Also called the ``moral" support for $h$ and
$\hat{h}$, respectively.} where ``most" of the information carried
by the function is located. As one may notice this is an {\it
absolute} way of quantifying the object since it relies on studying
the distribution of the $L^{\infty}$-norm of $h$ (and respectively
$\hat{h}$) and not on how $h$ may relate (interact) with some other
objects (functions) living in a given environment.

Reasoning in this spirit, (for $\f$ defined as above) we have that
the ``moral" support of $Q_b\f$ is given by the support of $\f$
(here we rely on the equality $\operatorname{supp}\:
Q_b\f=\operatorname{supp}\: \f$) \footnote{Remark that $Q_b$ is a
multiplication operator and hence preserves the time localization of
the object on which it acts - this being the main reason for which
we will split our operator $T$ in pieces that are compactly
supported in time.} while, with the notations from Figure 4, we
have that the ``moral support" of $\widehat{Q_b\f}$ is identified
with the frequency-interval $U$. At this point, we observe that we
lose the (global) Heisenberg principle\footnote{Throughout this
paper, we use the term ``Heisenberg principle" to refer to the
optimal Heisenberg localization, \textit{i.e.} the product of the
sizes of the time and frequency moral supports are comparable with
1.}, this being one of the main difficulties that was standing
against solving this conjecture.

One may improve this time-frequency portrait if one further
decomposes $\f$ in pieces which are better adapted to the
oscillation of the quadratic factor imposed by $Q_b$; more exactly,
writing
$$\f=\sum_{j}\f_j$$
with each $\f_j\in C_{0}^{\infty}$ and
$|\operatorname{supp}\f_j|\approx
\min((2b)^{-\frac{1}{2}},|\operatorname{supp}\f|)$ we squeeze the
previous localization to a sequence of area-one blocks concentrated
near the diagonal of the initial ``big"
rectangle. Now, even though on each such block - reflecting the
time-frequency portrait of a $\f_j$ - we recover the Heisenberg
principle, the parallelogram formed by their union still offers a
poor (global) localization of $Q_b\f$.  Using this viewpoint, one
cannot do better.

({\it q2}) The {\it relational} (relative) perspective, as the name
suggests, focuses on determining a contextual representation of our
object depending on how it interacts with other objects ``living" in
a given environment.

More exactly, in our case the environment is formed by
objects\footnote{Eventually, we will increase the complexity of these
objects by composing the symmetries.} like $M_{c}\f,\:Q_b\f$ and the
interaction is given by the scalar product in $L^{2}(\TT)$.

Now taking, for example, the interaction \beq\label{intM}
|\left\langle M_c\f ,M_{c'}\f \right\rangle | \eeq (here
$c,c'\in\R$) we see that, applying the (non-)stationary phase
principle, \eqref{intM} is controlled by a quantity depending on the
ratio of $|V|^{-1}$ and the distance between the lines $y=c$ and
$y=c'$ (where $\textrm{supp}\:\f=V$ and $\f$ is adapted to $V$). By varying $c'$, this
quantity suggests that (on the frequency side) the information
carried by $M_c\f$ should be localized ``near the line" $y=c$ and
that this information is roughly constant on intervals of length
$|V|^{-1}$. As a consequence we may interpret the relative
time-frequency localization of $M_c\f$ as being given by the region
(rectangle) centered near the line $y=c'$ of width  $|V|^{-1}$ (measured on
frequency axis) and with space support in the interval
$V$.\footnote{It is not surprising, in this case, that the relative
time-frequency picture coincides with the restrictive one described
above, given how the Fourier transform acts on modulation,
translation, and dilation.}

By analogy with the above description, we will now treat the
following interaction: \beq\label{intQ} |\left\langle Q_{b}\f
,Q_{b'}\f \right\rangle |\:. \eeq

As before, applying the (non-)stationary phase principle we remark
that \eqref{intQ} is controlled by a quantity depending on the ratio
of $|V|^{-1}$ and the distance\footnote{Here the appropriate notion of
distance is given
by $\sup_{x\in V}|2bx - 2b'x|$ rather than $\inf_{x\in V}|2bx - 2b'x|$.}
between the lines $y=2bx$
and $y=2b'x$ obtained by differentiating the polynomial phase.  This
fact invites us to think of the relative time-frequency localization
of $Q_b\f$ as being given by the region (parallelogram) centered
near the line $y=2bx$ of width $|V|^{-1}$ (measured on the frequency axis)
and with space support in the interval $V$.  Indeed, this
perspective will prove to give an accurate geometric representation
of the relations among our objects.\footnote{See Section 5.}

As a consequence, this should be the ``true"\footnote{Remark - see
Figure 4 - that using this approach, we recover (on each fiber) a
local Heisenberg principle.} time-frequency ``story" reflected in
pictures (see Figure 4); it is of {\it relative} nature since it
tells us about the interaction of $Q_b\f$ with an exterior object and
not about $Q_b\f$ itself.\footnote{For the remainder of the paper,
``time-frequency portrait" will refer to the \emph{relative}
representation described in (\textit{q2}).}  This time-frequency
interpretation can be regarded as a way of drawing pictures in which
besides the magnitude we also encode the oscillation of our
function.\footnote{The point is that while $|\widehat{Q_b\f}|$ is
big on the whole interval $U$, when tested against same-structure
functions (as in the expression $\left\langle \widehat{Q_b\f}
,\widehat{Q_{b'}\f} \right\rangle$) the oscillations of
$\widehat{Q_b\f}$ come into play, canceling out most of the
oscillations of $\widehat{Q_{b'}\f}$ up to the level given by the
interaction of the corresponding parallelograms (for further study
of this interaction behavior as well as for some other local
properties, see Section 5).}
\begin{figure}[!h]
\begin{center}
\epsfig{file=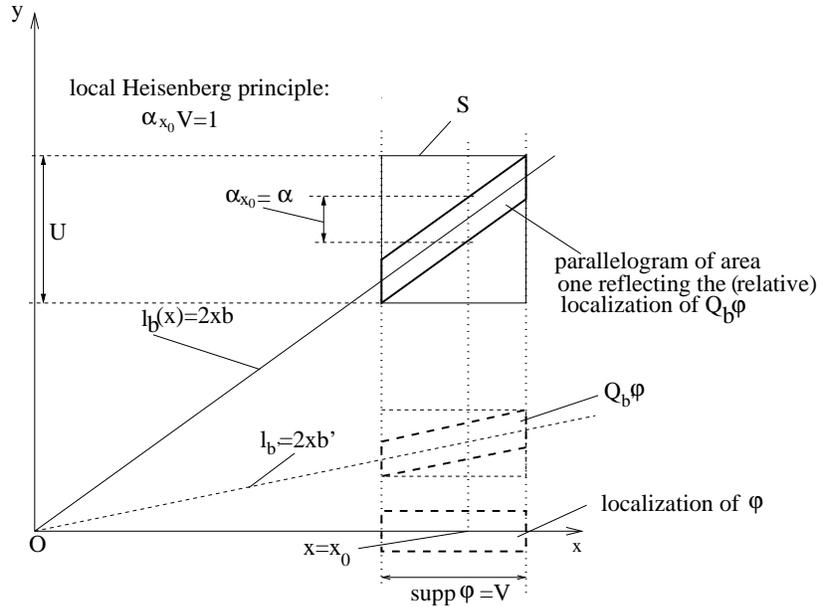,height=8cm}
\end{center} \caption{The (relative) time-frequency localization of $Q_b\f$}
\end{figure}\\

The moral of this story is that while $M_c$ translates the
time-frequency picture up and down, the operator $Q_b$ realizes a
shearing of the same picture.

The idea presented above will be essential in the proof of Theorem
1, and might be quite productive in a series of other problems
involving quadratic time-frequency analysis.

Now, coming back to our decomposition, if we let $M_c$ interfere
with $Q_b$ we obtain the ``elementary cell" of our operator modeled
in $M_cQ_b\f$; from the previous discussion, this will be considered
as being localized in a parallelogram of area one living near the
line $l(x)=c+2bx$ and with the same time localization as before (Figure
5).
\begin{figure}[!h]
\begin{center}
\epsfig{file=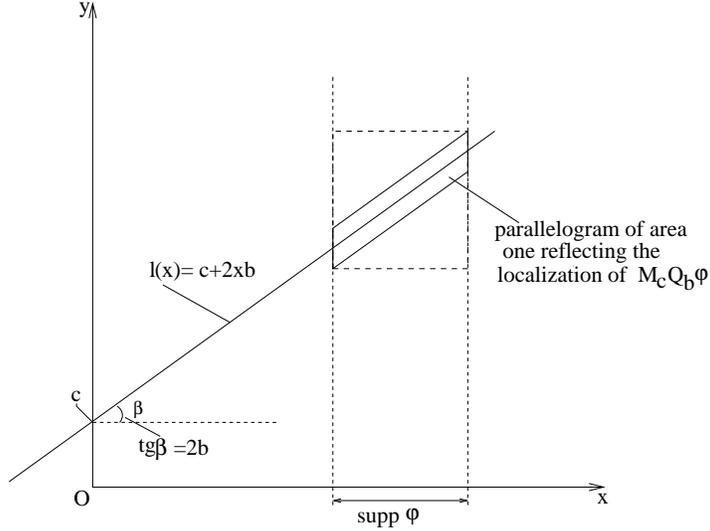,height=7cm}
\end{center} \caption{The (relative) time-frequency localization of $M_c Q_b\f$}
\end{figure}\\

Once we have gained this intuition, given the form \eqref{v6}, it is
natural to split $T$ in pieces that will be localized in the same
(relative) region as our ``elementary cells" $M_cQ_b\f$ mentioned
above. Consequently, we will divide our time-frequency plane in
parallelograms of area one as reflected in Figure 6.
\begin{figure}[!h]
\begin{center}
\epsfig{file=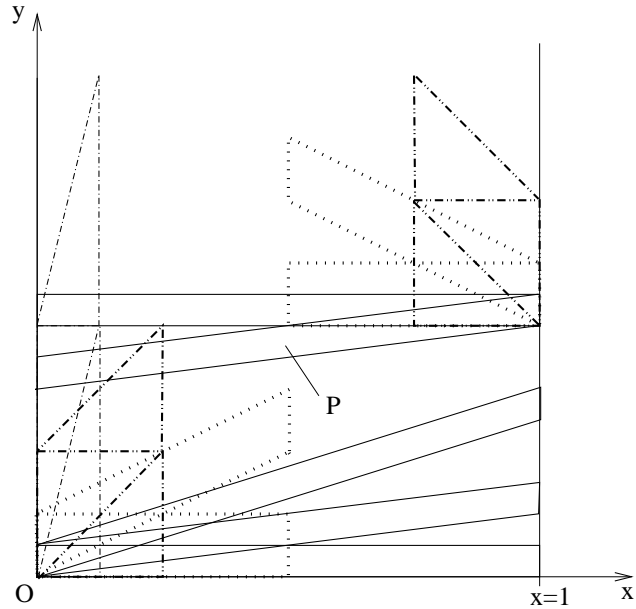,height=8cm}
\end{center} \caption{The time-frequency decomposition of the Quadratic Carleson Operator}
\end{figure}\\

The exact procedure will be described in Section 4, and will have as
a consequence
 $$T\:=\:\sum_{P\in\p}T_{p}\:,$$
 with each piece $T_p$ having the time-frequency picture
represented by the tile $P=[\a,\o,I]$ (see Section 3 for notations).

This way we have highlighted the dual nature of this problem: an
analytic formulation (providing $L^2$ bounds for a certain object)
visualized in terms of geometric interactions  of some families of
parallelograms (tiles). Consequently, there will be no surprise that
in the second stage of our program that we now initiate - the
selection algorithm - the geometric point of view in quantifying
different interactions among the ``small pieces" $T_P$ will play the
essential role.

Indeed, we start by defining a measurable map that assigns to each
point $x\in [0,1]$ a line $l_{x}\in \l$ in $\R^2$; then we can
regard $T_P f$ as  assigning the values (see Section 3 for
notations): \beq \label{tp1}
 \newline  x\:\stackrel{T_P f}{\longmapsto}\:0\:\:\:\:\:     \textrm{if} \:x\notin I\:\: \textrm{or}\:\: l_{x}\notin P
 \:,
 \eeq
 \beq \label{tp2}
 \newline  x\:\stackrel{T_P f}{\longmapsto}\: (\textrm{a quantity ``oscillating along}\:l_{x}")\:\:\:\:\:\textrm{if} \:x\in I\:\:\textrm{and}\:\:l_{x}\in P
 \:.
 \eeq

This way $T_P f$ (and similarly ${T_P}^{*} f$) encodes two different
types of information: $\newline$ - $\eqref{tp1}$ forces us to
consider the density of the ``flow" $\left\{l_{x}\right\}_{x\in I}$
through the tile $P$ (this concept will be made precise in Section 5
- see $\eqref{m}$ - and will be called the ``density factor" of P),
while $\eqref{tp2}$ implies that on Fourier side, the information
given by $\widehat{T_P f}$ is localized near the central line of $P$
denoted $l_P$. The interplay between these two features of $T_P$ (or
${T_P}^{*}$ ) will  be discussed in detail in Section 5, and it is
the key fact in providing good bounds for the expression
  \beq \label{Tp*}
 \left\|\sum_{P\in\p}{T_{P}}^{*}f\right\|_{2}^{2}=\sum_{P,P'\in\p}\left\langle {T_{P}}^{*}f ,
 {T_{P'}}^{*}f\right\rangle\:,
 \eeq
where here $\p$ is a certain finite collection of tiles and $f$ some
fixed element in $ L^2(\TT)$. In dealing with this problem, we first
need to understand the quantity
  \beq \label{tptp}
 \left|\left\langle {T_{P}}^{*}f , {T_{P'}}^{*}f\right\rangle
 \right|\:.
 \eeq

To obtain some intuition, we explain first the two possible
extreme cases:

- When $P=P'$ ({\it i.e.} the diagonal term) the relevant point of
view is given by $\eqref{tp1}$; this is natural since ${T_{P}}^{*}f$
and ${T_{P'}}^{*}f$ oscillate in the same region of the time
frequency-plane, making the information offered by $\eqref{tp2}$
useless. Consequently, the norm $\left\|T_P\right\|_2$ will measure
the density of $P$ (see $\eqref{tm}$ and $\eqref{m}$).

 - When $P$ and $P'$ are far apart from one another,
$\eqref{tptp}$ is small either due to the time localization of
${T_P}^{*}f$ or due to the relation $\eqref{tp2}$ that comes into
play by forcing $\widehat{{T_{P}}^{*}f}$ and
$\widehat{{T_{P'}}^{*}f}$ to have different ``moral supports".

Consequently, via $\eqref{tp1}$ and $\eqref{tp2}$ (which also
determines the time-frequency localization of $T_Pf$ and
${T_P}^{*}f$) we expect the following principle to be true:
  \beq\label{prin}
\begin{array}{rl}
&\operatorname{\: The\: magnitude\: of \:\eqref{tptp}\: is\: :}\\
&\operatorname{\:-\:big \:-}\: when\:P,\:P'\:\: have\:\:large\:\: overlaps\: \:and \:\: high\:\: density;\\
&\operatorname{\:-\:small\: - }\:when \:P,\:P'\:\: have\:\:small\:\:overlaps\:\:(are\:\: disjoint)\:\: or\\
 &low\:\: density.
 \end{array}
 \eeq

Now this principle  simultaneously  offers and demands a lot of
information:

(I) On the one hand, it suggests that to obtain good control of
$\eqref{Tp*}$ we may need to split the family of tiles $\p$ into
sub-collections $\p_j$ with each $\p_j$ having uniform
characteristics (all the tiles inside it must have comparable
densities and any interaction between two of them must have the same
degree of overlapping), estimate separately each
 $$\left\|{T^{\p_j}}^{*}f\right\|_2:=\left\|\sum_{P\in\p_j}{T_{P}}^{*}f\right\|_{2}$$
with bounds depending on the previously mentioned characteristics of
$\p_j$, and then sum them up for obtaining the desired global bound.

(II) On the other hand, it requires a clear formulation of the
concepts: $\newline$a) the density of a tile $\newline$ b) the
degree of the overlapping between two tiles .

Part (II) will be the object of our study in Section 5. While (II)
- a) will be straightforward, for (II) - b) we will introduce two
ways of measuring the corresponding concept: a \emph{qualitative}
one, by defining an ``almost" order relation between tiles - ``$\leq
$" - (Definition 3) and a \emph{quantitative} one, the actual
measurement of how much two tiles $P_1 , P_2$ overlap, that can be
recovered from the geometric factor of the pair $(P_1 , P_2)$
(Definition 1).

Now, guided by the observation made in (I), our proof will be based
on two propositions corresponding to the two main (geometric)
possibilities appearing in the study of a family of tiles (having
uniform density): Proposition 1 will treat the case where our family
consists of ``disjoint" ({\it i.e.} not comparable under ``$\leq$")
tiles, while Proposition 2 will deal with a family - called a
``forest" - that can be organized into a controlled number of
clustered sets of tiles ({\it i.e.} trees).

With this done we will proceed (roughly) as follows:

We will decompose $\P$ into $\bigcup_{n=0}^{\infty}\p_{n}$ with
$$\p_{n}=\left\{P\in\P\:|\:2^{-n-1}< \operatorname{the\:density\:factor\:of\:}P\leq 2^{-n}\right\}\:.$$

Using a combinatorial argument, we will further prove that $\p_n$
may be written as a disjoint union of at most $n$ sets,
$\bigcup_{j=1}^{n}\p_{nj}$, such that
$$\p_{nj}=\A_{nj}\:\cup\:\B_{nj}$$ where, for each $j$, $\A_{nj}$ is a family of
at most $n$ disjoint tiles  and $\B_{nj}$ is a forest. Now, denoting
$$T^{\p_{nj}}:=\sum_{P\in\p_{nj}}T_{P}\:,$$
Propositions 1 and 2 will imply that\footnote{Throughout this paper
we will denote with $\left\|T\:\right\|_2$ the operator norm of $T$
acting from $L^2$ to $L^2$.}
$$\left\|T^{\p_{nj}}\right\|_2\lesssim 2^{-n\eta}$$
for some absolute constant $\eta>0$, from which we conclude that
$$\left\|T\:\right\|_2\:=\left\|\:\sum_{n=0}^{\infty}T^{\p_n}\right\|_2\leq\sum_{n=0}^{\infty}\sum_{j=1}^{n}\left\|T^{\p_{nj}}\right\|_2\lesssim 1\:.$$
$\newline$

\section{\bf Notations}

Take the canonical dyadic grid in $[0,1]=\TT$ \footnote{For
convenience, from now on we may choose to identify $\TT$ with any
unit interval (not necessarily $[-\frac{1}{2},\frac{1}{2}]$). } and
in $\R$. A tile $P=[\a,\o,I]$ consists of dyadic (half open)
intervals $\a,\o \subset{\R}$ and $I \subset{[0,1]}$ with the
property that $\left|\a\right|=\left|\o\right|=\left|I\right|^{-1}$
(here $|I|=m(I)$ where $dm=dx$ stands for the Lebesgue measure on
$[0,1]$). The collection of all tiles $P$ \footnote{ For the
simplicity of notations, $P$ will encode two meanings (depending on
the context): a triple of intervals as defined above or the
parallelogram formed by these intervals in the time-frequency
plane.} will be denoted by $\P$.
\begin{figure}[!h]
\begin{center}
\epsfig{file=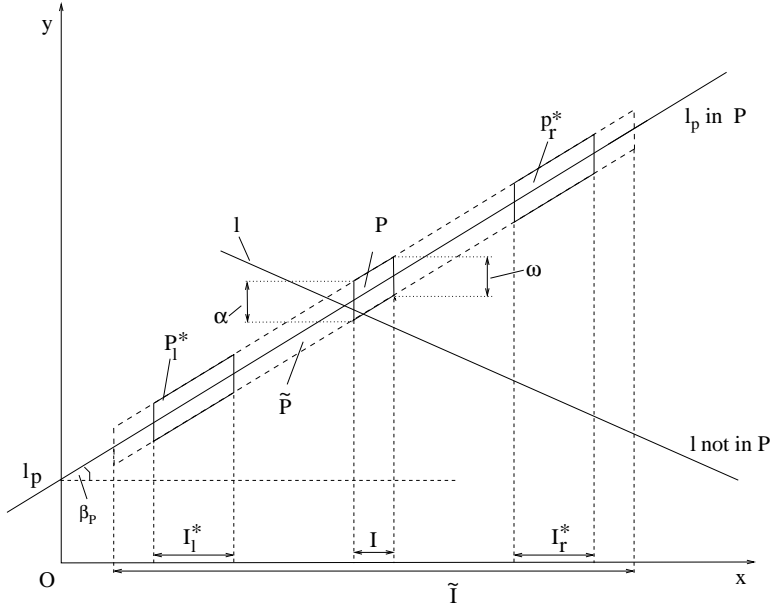,height=8cm}
\end{center} \caption{Notations}
\end{figure}\\
\indent If $I$ is any (dyadic) interval we denote by $c(I)$ the
center of $I$. Let $I_r$ be the ``right brother" of I, with
$c(I_r)=c(I)+|I|$ and $|I_r|=|I|$; similarly, the ``left brother" of
$I$ will be denoted $I_l$ with $c(I_l)=c(I)-|I|$ and $|I_l|=|I|$. If
$a>0$ is some real number, by $aI$ we mean the interval with the
same center $c(I)$ and with length $|aI|=a|I|$. Now for
$P=[\a,\o,I]\in\P$, we set $aP:=[a\a,a\o,I]$. Also, if $\p\subseteq
\P$ then by convention $a\p:=\left\{aP\:|\:P\in\p\right\}$.

Set $\l := \left\{ \textrm{all lines in the plane not parallel with
the $y$-axis} \right\}.$

\noindent Then, for each $P=[\a,\o,I]\in\P$ and $l\in\l$, we write
``$l\in P$" iff $l$ intersects both edges of $P$ which are parallel
with the $y$-axis. Also, for any tile $P$ as before, we will
associate the ``central line" $l_P$ - the unique line $l\in\l$ that
passes through the midpoints of the vertical edges (this line will
give the ``angle of $P$"- denoted $\beta_P$ and given by the formula
$\tan(\beta_P)=\frac{c(\o)-c(\a)}{|I|}$).

 Now, for $\beta\in\arctan(\Z)$, define
$$\p(k,\beta)=\left\{P=[\a,\o,I]\in\P\:|\:|I|=2^{-k}\:\:\:\&\:\:\:\beta_P=\beta\right\}\:.$$
Then this collection of disjoint (similar) parallelograms (tiles)
defines a partition of the band $\TT\times \R$. Fixing
$P=[\a,\o,I]\in\p(k,\beta)$, denote the ``upper brother" of $P$ by
$P_{u}=[\a_r,\o_r,I]\in\p(k,\beta)$; similarly, the lower brother of
$P$ will be $P_{l}=[\a_l,\o_l,I]\in\p(k,\beta)$.

 For any dyadic interval $I\subseteq [0,1]$ define the (non-dyadic) intervals
 $$I^{*}_{r}=[c(I)+\frac{7}{2}|I|,c(I)+\frac{11}{2}|I|)\:\:\:\&\:\:\:  I^{*}_l=[c(I)-\frac{11}{2}|I|,c(I)-\frac{7}{2}|I|)$$
 $I^{*}=I^{*}_{r}\cup I^{*}_{l}$  and  $\tilde{I}=13I$.

Similarly, for $P=[\a,\o,I]\in\P$ we adopt the following notations:
$\newline$ -$P^{*}_{r}$ for the tile (parallelogram of area two)
with time interval $I^{*}_{r}$ and the same central line $l_P$ as
$P$ $\newline$ -$P^{*}_{l}$ for the tile (parallelogram of area two)
with time interval $I^{*}_{l}$ and the same central line $l_P$ as
$P$.

The same procedure applies to $P^{*}$ and $\tilde{P}$ (see Figure
7).

Throughout the paper, for $f\in L^2(\TT)$, we denote by
$$Mf(x)=\sup_{x\in I}\frac{1}{|I|}\int_{I}|f|$$ the Hardy-Littlewood
maximal function associated to $f$.

If $\left\{I_j\right\}$ is a collection of pairwise disjoint
intervals in $[0,1]$ and $\left\{E_j\right\}$ a collection of sets
such that for a fixed $\d\in (0,1)$ \beq \label{sets}
 E_j\subset I_j\:\:\:\:\:\&\:\:\:\:\frac{|E_j|}{|I_j|}\leq\d\:\:\:\:\:\forall\:\:j\in
 \N\:,
\eeq then we denote

\beq \label{fmax} M_{\d}f(x):=\left\{
                        \begin{array}{rl}
                        \sup_{I\supset I_j}\frac{1}{|I|}\int_{I}|f|, \  \mbox{if} \  x\in E_j \\
                        0 \qquad, \  \mbox{if} \ x \notin E_j
                        \end{array} \right.
\:.\eeq

Remark that $\forall\:\:r>1$ we have

 \beq\label{v8}
 \left\|M_{\d}f\right\|_{r}^{r}\lesssim \d \left\|f\right\|_{r}^{r}
 \:.
 \eeq

For $A,\:B>0 $ we say $A\lesssim B\:(\gtrsim)$ if there exist an
absolute constant $C>0$ such that $A<CB\:(>)$; if the constant $C$
depends on some quantity $\d>0$ then we may choose to stress this
fact by writing $A\lesssim_{\d}B$.

If $C^{-1}A<B<CB$ for $C$ some small (positive) absolute constant
then we write $A\approx B$. For $x\in \R$ we set $\left\lceil
x\right\rceil:=\frac{1}{1+|x|}$.

The exponents $\eta$ and $\ep $  may change throughout the paper.

\section{\bf Discretization}

 Our aim is to ``properly" decompose the operator
 $$ Tf(x)=\sup_{b,c\in \R}|M_{c}Q_{b}H{Q_{b}}^{*}{M_{c}}^{*}f(x)|=\sup_{l\in\l}|T_{l}f(x)|$$
 where
 $$T_{l}f(x)=\int_{\TT}{\frac{1}{y}\:e^{i(l(x)y-by^2)}f(x-y)dy}$$ with
 $l\in\l$ given by $l(x)=c+2bx$.

Now linearizing\footnote{This procedure is often referred to as the
Kolmogorov linearization method.} T we can write
  $$Tf(x)=T_{l_x}f(x)=\int_{\TT}{\frac{1}{y}\:e^{i(l_x(x)y-b(x)y^2)}f(x-y)dy}$$
where by $l_x$ we understand a line in $\l$ given by
$l_x(z)=c(x)+2zb(x)$ where $c(\cdot)$ and $b(\cdot)$ are certain
measurable functions.

We start our decomposition by choosing $\psi$ to be an odd
$C^{\infty}$ function such that $\operatorname{supp}\:\psi\subseteq
\left\{y\in \R\:|\:2<|y|<8\right\}$ and having the property
$$\frac{1}{y}=\sum_{k\geq 0} \psi_k(y)\:\:\:\:\:\:\:\:\:\forall\:\:0<|y|<1\:,$$
where by definition $\psi_k(y):=2^{k}\psi(2^{k}y)$ (with $k\in \N$).
As a consequence, we deduce that
$$Tf(x)=\sum_{k\geq 0}T_{k}f(x):=\sum_{k\geq 0}\int_{\TT}e^{i\left\{l_{x}(x)y-b(x)y^2\right\}}\psi_{k}(y)f(x-y)dy\:.$$

Now for each $P=[\a,\o,I]\in\P$ let $E(P):=\left\{x\in I\:|\:l_x\in
P\right\}$. Also, if $|I|=2^{-k}$ ($k\geq0$) we define the operators
$ T_P$ on $L^2(\TT)$ by
$$T_{P}f(x)=\left\{\int_{\TT}e^{i\left\{l_{x}(x)y-b(x)y^2\right\}}\psi_{k}(y)f(x-y)dy\right\}\chi_{E(P)}(x)\:.$$

Clearly, as $P$ runs through
$\P_k:=\left\{P=[\a,\o,I]\in\P\:|\:|I|=2^{-k}\right\}$, for fixed
$k$, the $\left\{E(P)\right\}$ form a partition of $[0,1]$, and so
$$T_{k}f(x)=\sum_{P\in\P_{k}}T_{P}f(x)\:.$$

Consequently, we have
$$Tf(x)=\sum_{k\geq 0}T_{k}f(x)=\sum_{P\in\P}T_{P}f(x)\:.$$

This ends our decomposition.

We finish this section with several remarks.
 $\newline$1) Because we
want better separation properties between the support of $T_{P} f$
and that of $T^{*}_{P} f$ (for fixed $P$ and $f$), by further
splitting\footnote{We use here a partition of unity.} $\psi$ as:
$$\psi(y)=\sum_{j=1}^{13}\psi^{j}(y)$$
(with each $\psi^{j}$ an odd, smooth function with
$\operatorname{supp}\:\psi^{j}\subset\left\{1+\frac{j}{2}<\left|y\right|<
2+\frac{j}{2}\right\}$) we may assume (relabeling for example
$\psi^{6}$ with $\psi$) that
$$\operatorname{supp}\:\psi\subseteq \left\{y\in\R\:|\:4<|y|<5\right\}\:.$$

Consequently, for a tile $P=[\a,\o,I]$, the associated operator has
the properties
$$\operatorname{supp}\:T_P\subseteq I\:\:\:\:\:\:\:\:\:\:\:\&\:\:\:\:\:\:\operatorname{supp}\:T_P^{*}\subseteq \left\{x\:|\:3|I|\leq dist(x,I)\leq 5|I|\right\}=I^{*}$$
where here $T_P^{*}$ denotes, as usual, the adjoint of $T_P$.
$\newline$ 2) In what follows, (splitting
$\P=\bigcup_{j=0}^{9}\bigcup_{k\geq0}\P_{10k+j}$) we can suppose
that our collection $\P$ is sparse enough; namely, if
$P_j=[\a_j,\o_j,I_j]\in\P\:$with$\:j\in\left\{1,2\right\}$ such that
$|I_1|\not=|I_2|\:\:\:\operatorname{then}\:\:|I_1|\leq
2^{-10}|I_2|\:\:\:$or$\:\:\:|I_2|\leq 2^{-10}|I_1|$. $\newline$

\section{\bf Quantifying the interactions between tiles}

Our aim in this section is to isolate the appropriate quantities
that arise in the behavior of the expression
 \beq\label{IBT}
\left|\left\langle T_{P_1}^{*}f,T_{P_2}^{*}g\right\rangle\right|
\eeq and further to show how they control this interaction.

We begin our study by presenting a summary of the main properties
shared by the operator(s) involved in our considerations.

 $\newline$
 {\bf 5.1. Properties of $T_{P}$ and $T_{P}^{*}$}
$\newline$

 For $P=[\a,\o,I]\in\P$ with $|I|=2^{-k},\:k\in
\N$, we have \beq \label{v9}
\begin{array}{rl}
        &T_{P}f(x)=\left\{\int_{\TT}e^{i(l_{x}(x)y-b(x)y^2)}\psi_{k}(y)f(x-y)dy\right\}\chi_{E(P)}(x)\:,  \\
    &T_{P}^{*}f(x)=-\left\{\int_{\TT}e^{i(l_{x-y}(x-y)y+b(x-y)y^2)}\psi_{k}(y)\left(\chi_{E(P)}f\right)(x-y)dy\right\}\:.
\end{array}
\eeq

Notice that based on the previous interpretation of the symmetry
$Q_b$ (see Section 2), we may conclude:
\beq\label{loc}\begin{array}{rl} &\textrm{- the time-frequency
localization of $T_{P}$ is ``morally" given by the tile $P$;}\\
&\textrm{- the time-frequency localization of $T_{P}^{*}$ is
``morally" given by the (bi)tile $P^{*}$.}
\end{array}\eeq

Also, we have the pointwise estimate \beq\label{est}
\begin{array}{rl}
|T_{P}f(x)|\lesssim\frac{\int_{I^{*}}|f(y)|dy}{|I^{*}|}\chi_{E(P)}(x)\:
\end{array}\eeq
and the norm-estimate \beq\label{tm}
\left\|T_{P}\right\|_{2}\approx\left(\frac{|E(P)|}{|I|}\right)^{1/2}\:.
\eeq

 $\\\newline$ {\bf 5.2. Factors associated to a tile} $\newline$

Now, once we have understood what the main features of $T_{P}$ and
$T_{P}^{*}$ are, we will relate them to concepts regarding the
associated tile $P$. Indeed, taking into account relations
\eqref{est} and \eqref{tm}, and respectively \eqref{loc}, for a tile
$P=[\a,\o,I]$ we are naturally led to the following two quantities:

$\newline$a)$\:\:\:\:\:$     an {\it absolute} one (which may be
regarded as a self-interaction) that measures how many lines from
$\left\{\:l(x)\:\right\}_{x\in I}$ pass through $P$ relative to the
length of $I$; more exactly,
 we define the {\bf density (analytic) factor of $P$} to be the expression
\beq\label{m}
 A_{0}(P):=\frac{|E(P)|}{|I|}\:.
 \eeq

Notice from \eqref{tm} that $A_{0}(P)$ determines the $L^{2}$
norm of $T_{P}$. Consequently, we expect this quantity
to play an important role in organizing and estimating the family
$\left\{T_P\right\}_{P\in\p} $.

$\newline$ b)$\:\:\:\:\:\:$ a {\it relative} one (interaction of $P$
with something exterior to it) which is of geometric type: let be
$l\in \l$ a line and $P\in\P$ a tile as in Figure 8.
\begin{figure}[!h]
\begin{center}
\epsfig{file=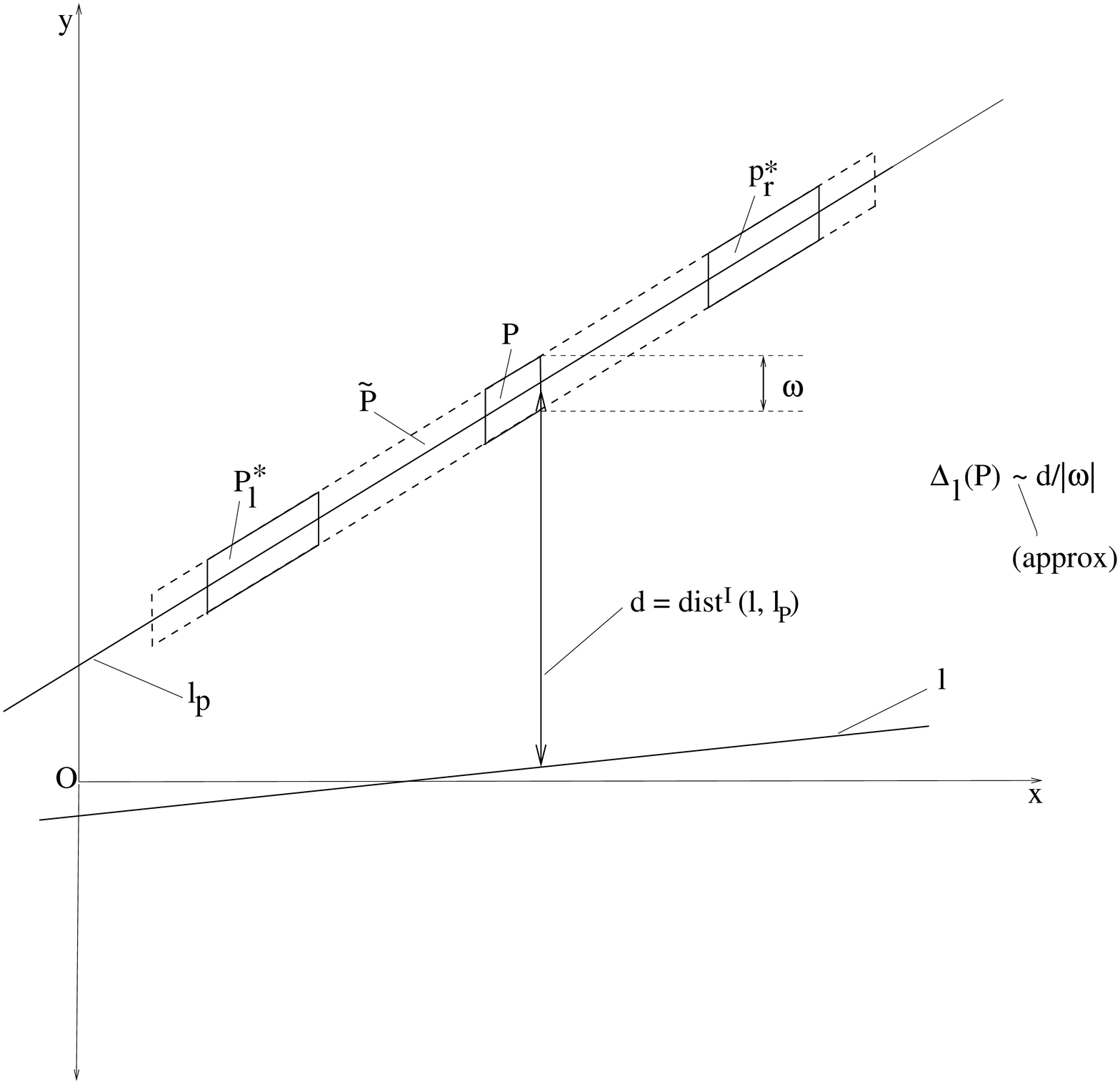,height=10cm}
\end{center} \caption{The geometric factor of $P$ with respect to $l$}
\end{figure}\\

For $l_{1},l_{2}\in\l$ we introduce the following notations:
 $$\operatorname{dist}_{x_0}(l_{1},l_{2})=\left|l_{1}(x_0)-l_{2}(x_0)\right|\:\:\:\:\&\:\:\:\:\operatorname{dist}^{A}(l_{1},l_{2})=\sup_{x\in A}\left\{\operatorname{dist}_{x}(l_{1},l_{2})\right\}\:.$$

Then we define the {\bf geometric factor of $P$ with respect to $l$}
to be the term

$$\left\lceil \Delta_l(P)\right\rceil\:,$$
where \beq\label{v1} \Delta_l(P):=\frac{\inf_{l_1\in
P}\left\{\operatorname{dist}^{I}(l,l_1)\right\}}{|\o|}\:. \eeq
$\newline$

{\bf 5.3. The resulting estimates} $\newline$

We now make the final step by observing how the above quantities
relate in controlling the interaction in \eqref{IBT}.

Given the heuristic \eqref{loc} and the form of \eqref{IBT}, we need
to quantify the relative position of $P_1^{*}$ with respect to
$P_2^{*}$.  To this end, we will need to adapt expression \eqref{v1}
to our context.\footnote{In the following we consider only the
nontrivial case $I_{P_1}^{*}\cap I_{P_2}^* \not=\es $.}

\begin{d0}\label{fact} Given two tiles $P_1$ and $P_2$ (suppose that $|I_1|\geq |I_2|$), we define the {\bf geometric factor of the pair ($P_1,P_2$)} by $$\left\lceil\Delta(P_1,P_2) \right\rceil\:,$$ where
$$\Delta(P_1,P_2)\:(=\Delta_{1,2}):=\frac{\inf_{{l_1\in P_1}\atop{{l_2\in P_2}}}\operatorname{dist}^{I_2}(l_{1},l_{2})}{|\o_2|}\:.$$
 \end{d0}

With these notations, remark that we have
$${\left\lceil\Delta_{1,2}\right\rceil}\approx\max\left\{{\left\lceil\Delta_{l_{P_1}}(P_2)\right\rceil},\:{\left\lceil\Delta_{l_{P_2}}(P_1)\right\rceil}\right\}\:.$$

$\newline$ \indent We will also need to define the ($\ep_0$-){\it
critical intersection interval} $I_{1,2}$ of the pair $(P_1,P_2)$ as
$$I_{1,2}=\left[x^{i}_{1,2}-\gamma_{1,2},x^{i}_{1,2}+\gamma_{1,2}\right]\cap
I_2^{*}\cap I_1^{*}$$ (see Figure 9).  Here $(x^i_{1,2} \, , \,
y^i_{1,2}):=l_{P_1}\cap l_{P_2}$ (if $l_{P_1}$ and $l_{P_2}$ are
parallel we set $x^i_{1,2} = \infty$), and $\gamma_{1,2}$ is chosen
to obey the relation \beq\label{gam}
\frac{\gamma_{1,2}}{\min(|I_1|,|I_2|)}=\left\lceil
\Delta_{1,2}\right\rceil^{\frac{1}{2}-\ep_0} \eeq for $\ep_0$  some
small fixed positive number.
\begin{figure}[!h]
\begin{center}
\epsfig{file=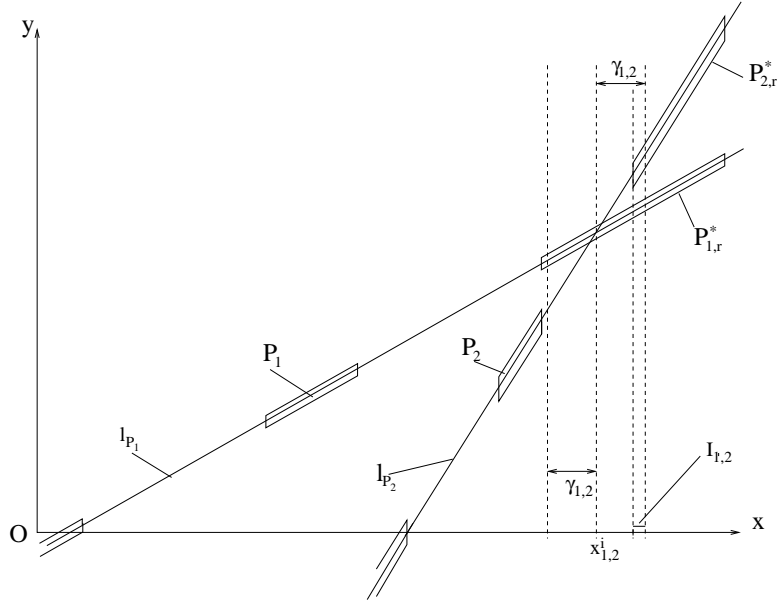,height=8cm}
\end{center} \caption{The interaction of two tiles}
\end{figure}\\
With these preparations done, we have the following result:

$\newline$ {\bf Lemma 0.} {\it Let be  $P_1\:,\:P_2\:\in\P$; then we
have \beq\label{v15} \left|\int
\tilde{\chi}_{I_{1,2}^c}T_{P_1}^{*}f\:\overline{T_{P_2}^{*}g}\:\right|\lesssim_{\:n,\:\ep_0}{\left\lceil
\Delta(P_1,P_2)\right\rceil}^{n}\:\frac{\int_{E(P_1)}
|f|\int_{E(P_2)}|g|}{\max\left(|I_1|,|I_2|\right)}\:\:\:\:\:\:\forall\:n\in
\N \eeq
 \beq\label{v16} \left|\int_{I_{1,2}}
T_{P_1}^{*}f\:\overline{T_{P_2}^{*}g}\:\right|\lesssim{\left\lceil
\Delta(P_1,P_2)\right\rceil}^{\frac{1}{2}-\ep_0}\:\frac{\int_{E(P_1)}
|f|\int_{E(P_2)}|g|}{\max\left(|I_1|,|I_2|\right)} \eeq where
$\tilde{\chi}_{I_{1,2}^c}$ is a smooth variant of the corresponding
cut-off.

Applying same methods for the limiting case $\ep_0=0$, we obtain
\beq\label{v17}
\left\|T_{P_1}{T}_{P_2}^{*}\right\|_2^{2}\lesssim\min\left\{\frac{|I_2|}{|I_1|},
\frac{|I_1|}{|I_2|} \right\}{\left\lceil
\Delta(P_1,P_2)\right\rceil}\:A_0(P_1)A_0(P_2)\:. \eeq}

The proof of Lemma 0 relies on the (non-)stationary phase principle
and is left to the reader.

\section{\bf The preparation - main ingredients}

As the title suggests, the role of this section is to present the
important concepts and results on which the proof of our theorem
relies.

We start on our way by introducing some quantitative and qualitative
notions that will help us later to organize our family of tiles.

The first step is to define a quantity that inherits relevant
features from both the analytic and geometric factors.

 \begin{d0}\label{mass}                     For $P=[\a,\o,I]\in\P$ we define the {\bf mass} of $P$ as being

\beq\label{v18} A(P):=\sup_{{P'=[\a',\o',I']\in\:\P}\atop{I\subseteq
I'}}\frac{|E(P')|}{|I'|}\:\left\lceil
\Delta(2P,\:2P')\right\rceil^{N} \eeq where $N$ is a fixed large
natural number.
\end{d0}

Next, we introduce a qualitative concept that characterizes the
overlapping relation between tiles.
 \begin{d0}\label{ord}
   Let $P_j=[\a_j,\o_j,I_j]\in\P$ with $j\in\left\{1,2\right\}$. We say that
$\newline$ - $P_1\leq P_2$       iff       $\:\:\:I_1\subseteq I_2$
and      $\exists\:\:l\in P_2\:\:s.t.\:\:l\in P_1\:,$ $\newline$ -
$P_1\trianglelefteq P_2$     iff     $\:\:\:I_1\subseteq I_2$  and
$\forall \:\:l\in P_2\:\:\Rightarrow\:\:l\in P_1\:.$
\end{d0}
\noindent{\bf Observation 1.} i) Remark that $\leq$ is not
transitive while $\trianglelefteq$ is. However, $\leq$ is not so far
from being a (partial) order relation; this may be encoded in the
fact that if $P_1\leq P_2$ then $2P_1\trianglelefteq 2P_2$.
$\newline$ii) Notice that the above definition can be meaningfully
extended (in the obvious manner) to any dilated tiles, {\it i.e.} it
makes sense to speak about $a_1 P_1\leq a_2 P_2$ and respectively
$a_1 P_1\trianglelefteq a_2 P_2$ (here $a_1,\:a_2>0$); in addition,
we say that $a_1 P_1\lneq a_2 P_2$ iff $a_1 P_1\leq a_2 P_2$ and
$|I_1|<|I_2|$. $\newline$iii) There is a nice connection between the
qualitative and quantitative concepts that measure the overlapping
of the tiles $P_1(\textrm{or}\:P_1^{*})$ and $P_2(\textrm{or}\:P_2^{*})$: if $I_1\subseteq I_2$ then
$\left\lceil
\Delta(P_1,P_2)\right\rceil=1\:\:\left(\Leftrightarrow\:\:\Delta(P_1,P_2)=0\right)\:\:\Leftrightarrow\:\:\left\{aP_1\leq
P_2\:\:\:\forall\:\:a>1\right\}$.

\noindent{\bf Observation 2.} Notice that the notion of mass of a
tile $P$ is dependent on the environment. This definition offers
many advantages, two of which we will mention here:
\begin{itemize}
\item the {\it monotonicity property (mp)}: if $P\leq P'$ (or $2P\trianglelefteq 2P'$) then $A(P)\geq A(P')$
\item the {\it smoothness property (sp)}: if $P$ and $P'$ are two tiles such that $I_P \approx I_{P'}$ ({\it i.e.} $2^{-a}I_P\subseteq I_{P'}\subseteq 2^{a}I_P$ for $a$
some small positive integer) and ${\left\lceil\Delta(P,P')\right\rceil}\approx 1$ then $A(P)\approx A(P')$.
\end{itemize}

\noindent{\bf Notation:} To avoid the boundary problems that may
arise from working with a single dyadic grid partition, we will
define the concept of the top (of a tree - see the next definition)
as being a set\footnote{This technicality is introduced only for
 smoothly handling the tree selection argument from Section 7.}
 of tiles $\tilde{P}=^{def}\{P^j\}_{j\in\{1,..s\}}$
 with $s\in \N$, $s\leq 4$ and $\{P^j=[\a^j,\o^j,I^j]\}_j$ having the
 properties:
 $\newline$1)$\:I^j=I^k\:\:\:\:\:\forall\:\:j,k\in\{1,\ldots s\}$
$\newline$2)$\:4P^j\leq 4P^{k}\:\:\:\:\forall\:\:j,k\in\{1,\ldots
s\}$ $\newline$ For $P\in\P$ we write $P\leq\tilde{P}$ iff
$\exists\:\:j\in\{1,..s\}$ such that $P\leq P^j$. In what follows,
it will also be convenient to work with a
representative\footnote{The reader may imagine a top as consisting
of only one (``fat") tile; indeed, in the following definitions and
results, the accent will always fall on a representative (of a top)
which may be regarded as a ``specialization" of the top itself.} of
the top $\tilde{P}$ - call it $P$ - which is some tile from
the collection $\{P^j\}_{j\in\{1,..s\}}$. $\newline$

Using the relation just defined, we now introduce the fundamental
(geometric) set-configuration that will govern most of our
reasonings.
 \begin{d0}\label{tree}

 We say that a set of tiles $\p\subset\P$ is a tree (relative to $``\leq"$) with top $\tilde{P}_0$ if
the following conditions are satisfied: $\newline
1)\:\:\:\:\:\forall\:\:P\in\p\:\:\:\Rightarrow\:\:\:\:\frac{3}{2}P\leq
\tilde{P}_0 $ $\newline 2)\:\:\:\:\:$if $P\in\p$ and
$\frac{3}{2}P_{u}\leq \tilde{P}_{0}$ then $P_u\in\p$ (analogously
for $P_l$) $\newline 3)\:\:\:\:\:$if $P_1,\:P_2\: \in\p$ and
$P_1\leq P \leq P_2$ then $P\in\p$
\end{d0}

\noindent{\bf Observation 3.} a) While conditions 1) and 3)
(appearing in the above definition) have clear
corespondents\footnote{The only difference appearing here is the
factor $\frac{3}{2}$ in 1) which is used for overcoming the boundary
problems that will arise later - see Section 7.2.} in \cite{2}, the
second condition - added here - is the extra twist that offers our
trees the advantage of being ``centered"\footnote{The central line
of the top (representative) splits the time-frequency representation
of our tree in two ``halves".}.

b) Sometimes we may exclude the (tiles forming the) top of the tree
from the collection $\p$. Also, we say that a tree has (top)
frequency line $l$ if $l$ is the central line of one of the tiles
(representative) belonging to the top. $\newline$

In this framework, we can state the results that will be used for
proving our theorem; their proofs will be postponed until Section 8.

\begin{p1}\label{prop1}
There exists $\eta\in(0,1/2)\:$ s.t. if  $\:\p\:$ is any given
family of incomparable tiles (i.e. no two of them can be related
through ``$\leq$") with the property that
$$ A(P)\leq\d\:\:\:\:\:\forall\:\:\:P\in\p$$
then
$$\left\|T^{\p}\right\|_2\lesssim\d^{\eta}\:.$$
\end{p1}

\begin{p1}\label{prop2}
Let $\left\{\p_j\right\}_j$ be a family of trees with tops
$\{\tilde{P}_j\}_j$ and respective representatives
$\{P_j=[\a_j,\o_j,I_j]\}_j$. $\newline$ Suppose that
$\newline\:1)\:\:\:\:A(P)<\d\:\:\:\:\:\:\:\forall\:j,\:P\in\p_j\:.$
$\newline\:2)\:\:\:\forall\:\:
k\not=j\:\:\&\:\:\forall\:\:P\in\p_j\:\:\:\:\:\:\:\:2P\nleq
2\tilde{P}_k\:.$ $\newline\:$3)    No point of $[0,1]$ belongs to
more than $K\d^{-2}$ of the $I_j$.

Then there is an absolute constant $\eta\in(0,\frac{1}{2})$ and a
set $F\subset\TT$ with $|F|\lesssim\d^{50}K^{-1}$ such that
$\forall\:f\in L^2(\TT)$ we have
$$\left\|\sum_{j}T^{{\p}_j}f\right\|_{L^2(F^{c})}\lesssim \d^{\eta}\log{K}\left\|f\right\|_2\:.$$

(Remark: Any collection of tiles $\p$ that can be represented as
$\cup_{j}\p_j$ with the family $\left\{\p_j\right\}$ respecting the
conditions mentioned above will be called a ``forest".)
\end{p1}

\noindent{\bf Observation 4.} One may notice the similarity between
the above propositions and the corresponding statements in \cite{2}
(Lemma 2 and Main Lemma); this is not surprising since the ``only"
difference between the quadratic case and the linear case is that we
have to deal with slanted rectangles. While the proof of Proposition
1 is basically the same as in \cite{2}, for the second proposition
we will have to deal with the extra overlaps of our parallelograms.
$\newline$

\section{\bf Proof of ``pointwise convergence"}

We now present the proof of Theorem 1.

 $\newline$ {\bf 7.1. Organizing the family of tiles}$\newline$

We start by breaking up $\P$ into $\bigcup_{n=0}^{\infty}\p_{n}$
where
$$\p_{n}=\left\{P\in\P\:|\:2^{-n-1}< A(P)\leq 2^{-n}\right\}\:.$$

Thus we have
$$T\:=\:\sum_{n=0}^{\infty}T^{\p_n}\:.$$

Here is the plan of our proof: $\newline$STEP 1 (the remaining part
of Section 7.1) - We modify each $\p_n$ so that the resulting set
gains a certain structure: all the elements inside it have
comparable mass and are clustered near some ``well-arranged" maximal
elements. $\newline$STEP 2 (Section 7.2) - Taking advantage of the
above-mentioned structure, we further show that each such $\p_n$ may
be decomposed (up to a negligible - in the sense of Proposition 1 -
family of tiles) into a certain number of forests.$\newline$ STEP 3
(Section 7.3) - Using Proposition 2, we will combine the estimates
for each forest into an estimate for the operator $T^{\p_n}$, which
allows us to obtain the desired bound for $T$.

As announced, we start the first part of our program by modifying
(cutting) some parts of the set $\p_n$. For this, we first define
$\left\{\bar{P}_{k}\right\},\:\bar{P}_{k}=[\bar{\a}_k,\bar{\o}_k,\bar{I}_k]$
to be the set of maximal triples with respect to $``\leq" $ that
obey the relation $\frac{|E(P)|}{|I_P|}\geq 2^{-n-1}$. Set
$\p_n^{0}$ to be \beq \label{triang}
 \p_n^{0}=\left\{P\in\p_n\:|\:\:\exists\:k\in \N \:s.t.\:\:\:\:\:4P\triangleleft \bar{P}_{k}\right\}
 \eeq
 and define also
$$\C_n=\left\{P\in\p_n\:|\:\operatorname{there\:are\: no\: chains}\:P\lneq P_{1}\lneq\ldots\lneq P_{n}\:\&\:\left\{P_j\right\}_{j=1}^{n}\subseteq\p_n\:\right\}\:.$$

With these notations we claim that
$$\p_n\setminus\C_n\subseteq\p_n^{0}\:.$$
Indeed, if $P\in\p_n\setminus\C_n$ then there exists a family of
tiles $\left\{P_j\right\}_{j=1}^{n}\subseteq\p_n$ such that $P\lneq
P_{1}\lneq....\lneq P_{n}$. Now, since $2^{-n-1}< A(P_n)\leq
2^{-n}$, we deduce that $\exists\:P'=[\a',\o',I']$ with
$I_n\subseteq I'$ such that $\frac{|E(P')|}{|I'|}\geq 2^{-n-1}$ and
$\Delta(P_n,P')<2^{n/3}$. From the maximality condition, we have
that $\:\exists\:\:k\in N\:\:\:s.t.\:\:\:P'\leq \bar{P}_{k}$ and so
$\Delta(P_n,\bar{P}_{k})<2^{n/2}$. On the other hand, from the chain
condition, we deduce that $\Delta(P,P_n)<3/2$ and $|\o_{P}|\geq
2^{n}|\o_{P_n}|\geq 2^{n}|\o_{\bar {P}_k}|$. Consequently, we have
that
$$\Delta(P,\bar{P}_{k})<3/2 + 2^{-n/2}\:,$$ which implies that $4P\triangleleft \bar{P}_{k}$ as we wanted.

 Let $\D_n\subseteq \C_n$ be the set such that
$\p_n\setminus\D_n=\p_n^{0}$; then $\D_n$ (or, in general, any
subset of $\C_n$) contains no (ascending) chains of length $n+1$ and
so breaks up as a disjoint union of a most $n$ sets
$\D_{n1}\cup\D_{n2}\cup\ldots\cup\D_{nn}$ with no two tiles in the
same $\D_{nj}$ comparable. Consequently, from Proposition 1, we have
$$\exists\:\:\eta\in(0,1/2)\:s.t.\:\:\left\|T^{\D_{nj}}\right\|_2\lesssim 2^{-n\eta}\:\:\:\:\:\:\:\forall \:\:\:\:j\in\left\{1,\ldots n\right\}\:,$$
which applied to $\D_n$ translates into \beq \label{D}
\left\|T^{\D_{n}}\right\|_2\leq\sum_{j=1}^{n}\left\|T^{\D_{nj}}\right\|_2\lesssim\sum_{j=1}^{n}2^{-n\eta}\lesssim
2^{-n\eta}\:. \eeq

As a consequence, we can now erase the set $\D_n$ without affecting
our plan. The resulting structure of the collection $\p_n^{0}$ will
help us later to further split our collection into forests, but for
the moment we turn our attention towards the set
$\left\{\bar{P}_{k}\right\}$, with the intention of obtaining a
rough bound for the counting function $N$ (defined below) associated
to the intervals $\left\{\bar{I}_{k}\right\}$. For this we notice
that $\left\{E(\bar{P}_{k})\right\}$ are pairwise disjoint, which
implies that $\sum_{k}|E(\bar{P}_{k})|\leq 1$. Now, using the
definition of $\bar{P}_{k}$, we deduce
$$\left\|N\right\|_{1}=\sum_{k}|\bar{I}_k|\leq 2^{n+1}|E(\bar{P}_{k})|\leq 2^{n+1}\:\operatorname{where}\:\:N(x)=^{def}\sum_{k}\chi_{\bar{I}_k}(x)\:.$$

Therefore the set defined as
$$G_n=\left\{x\in\TT\:|\:x\:\operatorname{is\:contained\:in\:more\:than\:}2^{2n}K\:\operatorname{of\:the\:}|\bar{I}_{k}|
\right\}$$
 has measure $|G_n|\lesssim (2^n K)^{-1}$. Because we want some control on the geometry of $\left\{\bar{P}_{k}\right\}$,  we will use $G_n$ for deleting more tiles from $\p_n^{0}$; indeed, if
$$\p_n^{G}=\left\{P=[\a,\o,I]\in\p_n^{0}\:|\:I\nsubseteq G_n\right\}\:,$$ we have that
\beq \label{exc}
T^{\p_n^{G}}f(x)\:=\:T^{\p_n^{0}}f(x)\:\:\:\:\:\:\:\:\:\:\forall\:f\in
L^{2}(\TT)\:\:\&\:\:x\in G_n^{c}\:. \eeq

(Since we have good control on the measure of $G_n$, we will focus
on estimating $T^{\p_n^{0}}$ only on $G_n^c$.)

 We delete from $\left\{\bar{P}_{k}\right\}$ all $\bar{P}_k$
with $\bar{I}_k\subseteq G_n$. Then the resulting set $\p_n^{G}$ has
the following properties: $\newline 1)\:\:A(P)\leq
2^{-n}\:\:\:\:\:\:\:\:\forall\:\:P\in\p_n^{G}\:,$ $\newline
2)\:\:\forall\: P\in\p_n^{G}\:\:\Rightarrow\:\exists\:k\in N
\:st\:\:\:\:\:4P\trianglelefteq \bar{P}_{k}\:,$ $\newline 3)\:\:$No
$x\in\TT$ belongs to more than $K 2^{2n}$ of the $\bar{I}_k$'s.

$\newline${\bf 7.2. Decomposing into forests}$\newline$

Now we shall prove that $\p_n^{G}$ decomposes\footnote{Up to a
family of chains with length controlled by an absolute constant.} as
a disjoint union of at most $M=2n\log K $ forests
$\B_{n0}\cup\B_{n1}\cup\B_{n2}\cup...\cup\B_{nM}\:,$ where each
$\B_{nk}$ satisfies the hypotheses of Proposition 2. In order to
make the decomposition, we first define
$$B(P)=\#\left\{j\:|\:4P\trianglelefteq \bar{P}_j\right\}\:\:\:\:\:\:\:\forall\:\:P\in\p_n^{G}\:.$$
Clearly $1\leq B(P)\leq 2^{M}$. Now let's define the sets
$$\p_{nj}:=\left\{P\in\p_n^{G}\:|\:2^j\leq B(P)< 2^{j+1}\right\}\:\:\:\:\:\:\:\:\:\:\forall\:j\in \left\{0,..M\right\}\:.$$

To better understand their behavior, we develop the following
procedure: fix a family of tiles $\p_{nj}$ as defined before and
$\newline$ 1) select the tiles
$\left\{P^{r}\right\}_{r\in\left\{1,\ldots
s\right\}}\subseteq\p_{nj}$ with the property that $\:4P^{r}$ are
maximal\footnote{Here we use the following convention: let be $\D$ a
collection of tiles; $P$ is maximal (relative to $``\leq"$) in $\D$
iff $\forall\:\:P'\in \D$ such that $P\leq P'$ we also have $P'\leq
P$.} elements with respect to the relation $``\leq"$ inside the set
$4\p_{nj}$. $\newline$ 2) from the maximality, we have  that \beq
\label{max1} 4P^{l}\leq 4P^{k}\:\:\Rightarrow\:\:I_l=I_k\:, \eeq
\beq \label{max2}
\forall\:P\in\p_{nj}\:\:\:\exists\:\:\:P^{l}\:\:\:\textrm{s.t.}\:\:\:\:4P\leq
4P^{l}. \eeq $\newline$ 3) from the definition of $\p_{nj}$ we
deduce
$$\operatorname{if}\:P\in\p_{nj}\:\:\operatorname{s.t.}\:\:\exists\:\:k\not=l\:\:\:\operatorname{with}\:\:\:\left\{{4P\trianglelefteq4P^{l}}\atop{4P\trianglelefteq 4P^{k}} \right.\:,\:\operatorname{then}\:\:\:\left\{{4P^{k}\leq4P^{l}}\:\atop{4P^{l}\leq 4P^{k}} \right.
\:.$$ $\newline$ 4) define $$\A_{nj}:=\left\{P\in\p_{nj}\:|\:
\forall\:\:P^{l}\:\:\Rightarrow\:\:\frac{3}{2}P\nleqslant
P^{l}\right\}\cup$$
$$\left\{P\:|\:\exists\:l\:st\:|I_P|=|I_{P^l}|\:,\:\frac{3}{2}P\leq P^{l}\:\&\:P\not=P^k\:\forall\:k\right\}=\A_{nj}^{1}\cup\A_{nj}^2$$

and set $$\p_{nj}=\A_{nj}\:\cup\:\B_{nj}\:.$$

Now, we claim that $\newline$a) $\A_{nj}$ can be split into a
controlled number of sets containing no chains (with respect to the
relation $``\leq"$). $\newline$b) the collection $\B_{nj}$ defines a
forest (up to a negligible family of tiles).

 We start with the proof of a), by
supposing that we can find $P_1,\:P_2\in A_{nj}^{1}$ such that
$P_1\lneq P_2$; suppose also (see \eqref{max2}) that \beq
\label{ch1}
 4P_2\leq 4P^{l}
\eeq for some $l$.

Now from the definition of $A_{nj}^{1}$ we have that
$\frac{3}{2}P_1\nleqslant P^{l}$, but observing that \beq
\label{ch3}
 |\o_{P_1}|\geq 2^{10}|\o_{P_2}|
 \eeq
  we contradict relation \eqref{ch1}.
The fact that the remaining set $\A_{nj}^2$ contains no chains comes
trivially from the maximality of the tiles
$\left\{P^{r}\right\}_{r\in\left\{1,..s\right\}}$.

 For part b), we proceed as follows:
$\newline$We choose $k\in\left\{1,..s\right\}$ and define \beq
\label{tb}
 S_k=\left\{P\in\B_{nj}\:|\:\frac{3}{2}P\lneq P^{k}\right\}\:.
\eeq We now collect all $\left\{P^k\right\}_k$ for which $S_k=\es$
and erase them using Proposition 1. Consequently, by relabeling the
remaining maximal tiles we can always suppose that for each $P^k$ we
have $S_k\not=\es$ and that $\B_{nj}=\bigcup_{k}\left\{S_k\cup
P^k\right\}$. Further, we want to study the separation properties of
the family $\left\{S_k\right\}_k$. For this, we first introduce the
following relation: we say that
$$S_k\propto S_l$$
 if and only if $\exists\: P_1\in S_k$ and
$\exists\: P_2\in S_l$ such that $2P_1\leq 2P_2$ or $2P_2\leq 2P_1$.

 With this done, we first claim that \beq \label{rel}
S_k\propto
S_l\:\:\:\Rightarrow\:\:\:4P^k\leq4P^l\:\:\Rightarrow\:\:I^k=I^l\:.
\eeq

Indeed, suppose that $S_k\propto S_l$, and so (without loss of
generality) we know that $k\not=l$ and there are $ P_1\in S_k$ and
$P_2\in S_l$ such that $2P_1\leq 2P_2$. Then, since
$\frac{3}{2}P_2\leq P^l$ and $|\o_1|\geq |\o_2|\geq 2^{10} |\o^l|$,
we must have $4P_1\trianglelefteq 4P^l$. On the other hand, since $
P_1\in S_k$, we also have $4P_1\trianglelefteq 4P^k$, but this
 forces (see 3)) $4P^k\leq4P^l$.

We now construct the sets $$\bar{S}_k:=S_k\cup P^k\:\:\:k\geq1$$ and
observe that with a similar reasoning as in \eqref{rel} we obtain
\beq \label{eqrel} \bar{S}_k\propto
\bar{S}_l\:\:\:\Rightarrow\:\:\:4P^k\leq4P^l\:\:\Rightarrow\:\:I^k=I^l\:.
\eeq

The point is that with respect to $\left\{\bar{S}_k\right\}_k$,
$\propto$ becomes an equivalence relation. Indeed, let us check  the
transitivity of our relation. Suppose that $\bar{S}_k\propto
\bar{S}_l\propto \bar{S}_m$. Now, since $\bar{S}_k\propto
\bar{S}_l$, we deduce from \eqref{eqrel} that $4P^k\leq 4P^l$, and
since $S_k\not=\es$ we also have that $\exists\:P_1\in S_k$ with
$\frac{3}{2}P_1\lneq P^{k}$. On the other hand, from
$\bar{S}_l\propto \bar{S}_m$, we have that $4P^l\leq 4P^m$. Putting
these facts together, we have that $10P^k\leq 10P^m$, $I^k=I^m$, and
since $|\o_1|\geq 2^{10} |\o^k|$ we deduce $2P_1\trianglelefteq
2P^m$, which proves our claim.

Now let $\hat{k}:=\left\{m\:|\:\bar{S}_m\propto \bar{S}_k \right\}$
(observe that the size of the orbit of each $k\:(\bar{S}_k)$ is at
most 4). Denote
$$\hat{S}_k:=\bigcup_{m\in\hat{k}}\bar{S}_m\:.$$

Now, choosing a unique representative in each equivalence class, and
relabeling the resulting elements in a consecutive order, we deduce
that $\hat{S}_k\cap \hat{S}_l=\es $ for any $k\not= l$, which
implies $\left\{\hat{S}_k\right\}_{k}$ is a partition of $\B_{nj}$.

We need some final modifications to each set $\hat{S}_k$. First, we
denote by $\tilde{P}^k$ the set of all the maximal tiles
$\left\{P^l\right\}_l$ contained in the collection $\hat{S}_k$; now,
using Proposition 1, we delete, for each $k$, all the elements
belonging to $\tilde{P}^k$. Secondly, we define (for each k) the set
$$\hat{S}^{min}_k=\left\{P\in\hat{S}_k\:|\:\forall\:P'\in\hat{S}_k\:st\:I_{P'}\cap I_{P}\not=\es\:\Rightarrow\:I_{P}\subseteq I_{P'}\right\}\:.$$

By construction, the set $\left\{\hat{S}^{min}_k\right\}_k$ contains
only pairwise disjoint tiles, so, again applying Proposition 1, we
can erase the set $\hat{S}^{min}_k$ from each $\hat{S}_k$ and
consider $\B_{nj}=\bigcup_{k}\hat{S}_k$. In what follows we will
prove that each $\hat{S}_k$ is a tree with top $\tilde{P}^{k}$.

Indeed, fix a collection $\hat{S}_k$ of tiles; we will now verify
conditions 1)-3) in Definition 4. Take $P\in \hat{S}_k$; first
observe that 1) holds trivially since by construction
$\frac{3}{2}P\lneq \tilde{P}^{k}$. Suppose now that $P\in\hat{S}_k$
with $\frac{3}{2}P_u\leq \tilde{P}^{k}$. Then to show $P_u
\in\hat{S}_k$ it is enough to prove that $P_u\in \B_{nj}$. For this,
we need first to prove that $P_u\in \p_n^{G}$. Since
$\frac{3}{2}P_u\leq P^{k}$ (for some $P^k$ an element of
$\tilde{P}^k$) the above statement reduces to $P_u\in \p_n$. But we
know that $2P_u\trianglelefteq 2P^k$, and since $A(P^{k})>2^{-n-1}$,
using {\it (mp)}, we deduce that also $A(P_u)>2^{-n-1}$. At this
point, we recall that (following the previous procedure) \beq
\label{min} \exists\:\: P_0\in\hat{S}^{min}_k \:\:\:\:st\:\:\:\:
I_{P_0}\subsetneq I_P\:. \eeq

Consequently, using \eqref{min}, we have $2P_0\trianglelefteq 2P_u$,
and so by {\it (mp)} $A(P_u)\leq A(P_0)\leq 2^{-n}$. For the second
part, we need $P_u\in\p_{nj}$, but this comes from the fact that
 $$4P_u\trianglelefteq 4P^{k}\:\:\:\:\Rightarrow\:\:\:\:B(P_u)\geq 2^{j}$$
 and
 $$4P_0\trianglelefteq 4P_u\:\:\:\:\Rightarrow\:\:\:\:B(P_u)<
 2^{j+1}\:.$$

From this, we conclude that $P_u\in \hat{S}_k$, so 2) is true.

The convexity condition 3) is trivial since if $P_1<P_2<P_3$ with
$P_1\:\:\&\:\:P_3\in \hat{S}_k $ we have
$\frac{3}{2}P_1\trianglelefteq \frac{3}{2}P_2\trianglelefteq
\frac{3}{2}P_3\leq\frac{3}{2}\tilde{P}^k$, which implies $P_2\in
\hat{S}_k$.

 Consequently, we have proven $\hat{S}_k$ is a tree with
top $\tilde{P}^k$. Now, from the previous considerations, we have
that
$$\hat{S}_k\propto\hat{S}_l\:\:\Rightarrow\:\:k=l\:,$$
and since $\B_{nj}=\bigcup_{k}\hat{S}_k$ we deduce that $\B_{nj}$
becomes a forest as defined in Proposition 2.

$\newline${\bf 7.3. Ending the proof} $\newline$

 Now, we may conclude as in \cite{2}.

We first apply Proposition 2 for each family $\B_{nj}$ and obtain
that
$$\left\|T^{{\p}_{nj}}f\right\|_{L^{2}(F^{c}_{nj})}\lesssim 2^{-n\eta}\log K \left\|f\right\|_2\:,$$
where $F_{nj}$ is a small set with measure $|F_{nj}|\lesssim
{2^{n}K}^{-1}$. As a result, denoting $F_{n}=\cup_{j} F_{nj}$, we
have that \beq \label{e}
\left\|T^{{\p}_n^{G}}f\right\|_{L^{2}(F^{c}_{n})}\leq
\sum_{j=1}^{2n\log
K}\left\|T^{{\p}_{nj}}f\right\|_{L^{2}(F^{c}_{nj})}\lesssim n
2^{-n\eta}{(\log K)}^{2}\left\|f\right\|_2 \eeq with $|F|\lesssim
\frac{n\log K}{2^n K}$.

Therefore, combining  \eqref{D}, \eqref{exc} and \eqref{e}, we
deduce
$$\left\|T^{{\p}_n}f\right\|_{L^{2}(E^{c}_{n})}\lesssim n 2^{-n\eta}{(\log K)}^{2}\left\|f\right\|_2\:,$$
where $E_{n}=F_{n}\cup G_{n}$ still has measure $\lesssim
\frac{n\log K}{2^n K}$.

 Summing now over $n$, we obtain
$$\left\|Tf\right\|_{L^{2}(E^{c})}\lesssim {(\log K)}^{2}\left\|f\right\|_2$$
with $E=\cup_{n}E_n$ and $|E|\lesssim\frac{\log K}{K}$.

In conclusion, given $\g\:>\:0$, we have that for all $K>100$
$$\left|\left\{|Tf(x)|\:>\:\g\right\}\right|\leq\frac{\left\|Tf\right\|^{2}_{L^{2}(E^{c})}}{{\g}^2}\:+\:|E|\lesssim
{(\log K)}^{4}\frac{\left\|f\right\|_2^{2}}{{\g}^2}\:+\:\frac{\log
K}{K}\:.$$

Now, if we pick $K$ to minimize the right-hand side, we arrive at
the relation
$$\left|\left\{|Tf(x)|\:>\:\g\right\}\right|\lesssim_{\ep}\left(\frac{\left\|f\right\|_2}{{\g}}\right)^{2-\ep}\:\:\:\:\:\:\:\:\:\:\:\forall\:\:\:\:\ep\in\:(0,2)\:,$$
which further implies
$$\left\|Tf\right\|_p\lesssim_{p} \left\|f\right\|_2\:\:\:\:\:\:\:\:\:\:\:\:\:\forall\:\:\:\:p\:<\:2\:,$$
ending the proof of our theorem.

\section{\bf Some technicalities - the proofs of Propositions 1 and
2}

$\newline$\indent We now present the proofs of the
statements made in Section 6.

$\newline$ {\bf Proof of Proposition 1}
$\newline$

 The basic idea of our proof relies on combining the
$TT^{*}$ and maximal methods. Indeed, once we have expressed the
norm of our operator as a sum of interactions among ``small pieces"
$T_P $, we split it in two terms:

- for the first one (close to the diagonal) we use some maximal
methods since all our pieces $T_P $ ``oscillate" in the same region
of the time-frequency plane,

- for the second one (far from the diagonal) we take advantage of
the orthogonality of our terms, which is reflected in the smallness
of the resulting geometric factors.

$$\int_{\TT}\left|\left ({T^{\p}}\right )^{*}f(x)\right|^2dx\lesssim\left|\sum_{P^{\prime}\in\p\atop{P^\prime=[\a^\prime,\o^\prime,I^\prime]}}
\int_{\TT}f(x)\left\{\sum_{P=[\a,\o,I]\in\:\p\atop{\left|I\right|\leq\left|I^\prime\right|}}
\overline{T_{P^\prime }T^{*}_{P}f}(x)\right\}dx\:\right|$$
$$\lesssim\sum_{P^\prime\in\p}\int_{E(P^\prime)}|f|\left\{\sum_{P\in a(P^\prime)}{\left\lceil
{\Delta}(P,P^\prime)\right\rceil}^{1/2}\frac{\int_{E(P)}\left|f\right|}{\left|I^\prime\right|}\right\}$$
$$+\:\sum_{P^\prime\in\p}\int_{E(P^\prime)}|f|\left\{\sum_{P\in b(P^\prime)}{\left\lceil
{\Delta}(P,P^\prime)\right\rceil}^{1/2}\frac{\int_{E(P)}\left|f\right|}{\left|I^\prime\right|}\right\}=^{def}
A\:+\:B$$ where for the third inequality we used the estimate (cf.
Lemma 0)
$$\left|T_{P^\prime}T^{*}_{P}f(x)\right|\lesssim{\left\lceil {\Delta}(P,P^\prime)\right\rceil}^{1/2}
\frac{\int_{E(P)}\left|f\right|}{\left|I^\prime\right|}\chi_{E(P^\prime)}(x)$$
together with the following notations:
$$a(P^\prime)=\left\{P=[\a,\o,I]\in\p\:,\:\left|I\right|\leq\left|I^\prime\right|\:\:\&\:\:I^{*}\cap{I^\prime}^{*}\not={\es}\:\:|\:\:\Delta(P,P^\prime)\leq\d^{-2\ep} \right\}\:,$$
$$b(P^\prime)=\left\{P=[\a,\o,I]\in\p\:,\:\left|I\right|\leq\left|I^\prime\right|\:\:\&\:\:I^{*}\cap{I^\prime}^{*}\not={\es}\:\:|\:\:\Delta(P,P^\prime)\geq\d^{-2\ep} \right\}\:.$$
(Here $\ep\in (0,1)$ is some fixed constant.)

 Now the second term is easy to estimate:
$$ B\lesssim\sum_{P^\prime\in\p}\int_{E(P^\prime)}\left|f(x)\right|\left\{\frac{\d^{\ep}}{\left|I^\prime\right|}\sum_{P\in b(P^\prime)}\int_{E(P)}\left|f\right|\right\}dx\leq$$
$$\d^{\ep}\sum_{P^\prime\in\p}\int_{E(P^\prime)}\left|f(x)\right|\left\{\frac{1}{\left|I^\prime\right|}\int_{\tilde{I^\prime}}\left|f\right|\right\}dx
\leq\d^{\ep}\int_{\TT}\left|f(x)\right|Mf(x)dx\lesssim\d^{\ep}\int_{\TT}\left|f\right|^2\:.$$

For the first term we use the following Carleson measure-type estimate:
\beq \label{cm}
\sum_{P\in a(P^\prime)}|E(P)|\lesssim \d^{1-100\ep}\:|I^\prime|\:,\eeq
which is a consequence of the smoothness property {\it (sp)} of the
mass $A(P)$. Indeed, define
$$\J(P^\prime)=\{I\:|\:\exists\:P=[\a,\o,I]\in a(P^\prime)\:\}\:.$$
Let $\J_{min}(P^\prime)$ be the set of minimal (with respect to inclusion) intervals inside $\J(P^\prime)$, and define
\begin{align*}
\check{\J}(P^\prime):=\big\{I\subset 30 I^\prime\:|\: & \textrm{Exactly one of the left or right halves}\\
&\textrm{of $I$ contains an element of $\J_{min}(P^\prime)$} \big\} \: \cup \J_{min}(P^\prime)\:.
\end{align*}
Finally, set $$\breve{a}(P^\prime)=\left\{P=[\a,\o,I]\in\P\:\:|\:\:I\in\check{\J}(P^\prime)\:\:\&\:\: \Delta(P,P^\prime)\leq\d^{-2\ep} \right\}\:.$$
Then using the property {\it (sp)} \footnote{With an appropriate choice of $N$ in the definition of $A(P)$.} and the fact that any two tiles inside $\p$ are not comparable we have
$$\sum_{P\in a(P^\prime)}|E(P)|\leq\sum_{P\in \breve{a}(P^\prime)}|E(P)|\lesssim \d^{1-100\ep}\:\sum_{I\in \check{\J}(P^\prime)}|I|\lesssim \d^{1-100\ep}\:|I^\prime|\:,$$
which gives us the desired estimate (\ref{cm}).

Now set $E_{P^\prime}:=\cup_{P\in a(P^\prime)}E(P)$; using
H\"older's inequality for some fixed $1<r<2$, we deduce
$$A\lesssim\sum_{P^\prime\in\p}\frac{\int_{E(P^\prime)}\left|f\right|}{\left|I^\prime\right|}\int_{E_{P^\prime}}|f|
\lesssim\sum_{P^\prime\in\p}\int_{E(P^\prime)}\left|f\right|\left(
\frac{|E_{P^\prime}|}{|I^\prime|}\right)^{1-\frac{1}{r}}\left(\frac{\int_{I^\prime}|f|^{r}}{|I^\prime|}\right)^{\frac{1}{r}}$$
$$\lesssim\d^{1-1/r-100\ep}\sum_{P^\prime\in\p}\int_{E(P^\prime)}\left|f(y)\right|f^{*}_{r}(y)dy\:$$
$$\lesssim\d^{1-1/r-100\ep}\int_{\TT}[f^{*}_{r}(y)]^2dy\:\lesssim\:\d^{1-1/r-100\ep}\left\|f\right\|_2^{2}$$
(Here $f^{*}_{r}(x)=\sup_{x\in
I}(\frac{\int_{I}\left|f(y)\right|^{r}dy}{\left|I\right|})^{1/r}$
designates the Hardy-Littlewood maximal function of order $r$.)

The conclusion of our proposition now follows, by properly choosing
$\ep>0$.

\begin{flushright}
$\Box$
\end{flushright}

The remainder of the section will be dedicated to proving
Proposition 2. The natural approach is to obtain control on:

- the single tree estimate (Lemma 1),

- the interaction between (separated) trees (Lemmas 2 and
3).$\newline$

We now start the study of the other\footnote{As opposed to the
structure of the family of tiles appearing in Proposition 1.}
extremal geometric configuration, namely the tree. For Lemma 1, due
to the structure of our family, the geometric factors will play no
role, the entire effort being concentrated on properly using the
(uniform) density condition and the mean zero property.

\begin{l1}\label{tr}
Let $\delta>0$ be fixed and let $\p\subseteq\P$ be a tree with top
$\tilde{P}_0$, representative $P_0=[\a_0,\o_0,I_0]$, and frequency
line $l_0$ and such that
$$A(P)<\delta\:\:\:\:\forall\:\:\:P\in\p\:.$$ Then
$$\left\|T^\p\right\|_2\leq{\delta^{1/2}}\:.$$
\begin{proof}

The essence of the proof\footnote{Our case is a ``quadratic
perturbation" (that realizes a shearing) of the linear tree case
presented in \cite{2}.} below relies on the following outlook:
\beq\label{concept}
\begin{array}{rl}
&\textrm{``For $\p$ a tree, the associated operator $T^\p$ behaves}\\
&\textrm{like the (maximal) Hilbert transform."}
\end{array}
\eeq

[{\it Remark.} The easiest way to understand this heuristic is to
take a particular instance of $\p$ ($T^\p$): suppose that the top
$P_0$ stays on the real axis, that $l_x\equiv 0$ for any $ x\in I_0$
and that all the minimal tiles in the collection $\p$ (we may assume
$\p$ finite) are at the same scale.  Then, from the convexity
condition 3) in Definition 4, we remark that $\exists\: k_0,\:k_1\in
\Z$ such that, for $x\in I_0$, \beq\label{simplt}
T^{\p}f(x)=\sum_{k_0\leq k \leq k_1}\int{\psi_k(y)f(x-y)dy}, \eeq
{\it i.e.} $T^{\p}$ is a truncation of \eqref{hilb}.]

To make this precise, we will further show that \eqref{simplt} is
always true locally, on $supp\: T^{\p}$. Indeed, let $\newline
\begin{array}{rl}
&k_0(x)=\inf\left\{k\in \N\:|\:\exists\:P\in\p\:\:st\:\:|I_p|=2^{-k}\:\:\&\:\:\chi_{E(P)}(x)\not=0\right\}\:,\\
&k_1(x)=\sup\left\{k\in
\N\:|\:\exists\:P\in\p\:\:st\:\:|I_p|=2^{-k}\:\:\&\:\:\chi_{E(P)}(x)\not=0\right\}\:.
\end{array}$

Using the convexity condition we then deduce \beq\label{v20}
T^{\p}f(x)=\sum_{k_0(x)\leq k \leq
k_1(x)}\int{e^{i(l_x(x)y-b(x)y^2)}\psi_k(y)f(x-y)dy}\:. \eeq

Since we also want to obtain some decay, we need to take advantage
of the ``mass" of our tree. For this, the key fact is to observe
that heuristically our operator behaves as follows:
 \beq \label{heur}
 \begin{array}{rl}
&x \ \stackrel{T^{\p} f}{\longmapsto} \ 0 \ \ \textrm{if} \ x\in I^j\setminus E^j \ ,\\
&x \ \stackrel{T^{\p} f}{\longmapsto} \ \textrm{constant$(j)$} \ \
\textrm{if} \ x\in E^j \ ,
\end{array}
 \eeq
where here the sets $\left\{I^j\right\}_{j}$ and
$\left\{E^j\right\}_{j}$ obey the conditions\footnote{Indeed, one
can define $\left\{I^j\right\}_{j}$ to be the maximal dyadic
intervals contained in $I_0$ that satisfy
$$\frac{|E(l_0,I)|}{|I|}>100\d $$
 where $E(l_0,I):=\left\{x\in I\:|\:\operatorname{dist}^{I}(l_x,l_0)<2|I|^{-1}\right\}$.
Now setting $\bar{E}^j=E(l_0,\bar{I}^j)$ and $E^j=\bar{E}^j\cap I^j$
and making use of {\it (sp)} one concludes that if
$P=[\a,\o,I]\in\p$ with $I\cap I^j\not=\emptyset$ then
$\bar{I}^j\subseteq I\:\:\:\&\:\:\:E(P)\cap I^j\subseteq E^j$ (here
$\bar{I}^j$ is the dyadic interval containing $I^j$ and having the
length twice as big).}: $\newline$ $\:\:\:-\left\{I^j\right\}$ is a
partition of $I_0$\:, $\newline \:\:\:-E^j\subseteq I^j\:
\textrm{and}\:\frac{|E^j|}{|I^j|}\lesssim \d $\:.

Now, combining the views offered by \eqref{concept} and \eqref{heur}
we proceed as follows: To come closer to \eqref{simplt}, our first
step is to move our tree near the real axis: set
$\T^{\p}={Q_{b_0}}^{*}{M_{c_0}}^{*}T^{\p} M_{c_0}Q_{b_0}\:$ and
$g(x)=M_{c_0}^{*}Q_{b_0}^{*}f(x)$
 (here $l_0(z)=c_{0}+2b_{0}z$ is the central line of $P_0$).
Then, for $x\in I_0$ fixed, we have
$$\left|T^{\p}f(x)\right|=\left|\T^{\p}g(x)\right|\leq$$
$$\sum_{k_0(x)\leq k \leq k_1(x)}\left\{\int_{\TT}\left|e^{i\left\{(l_x(x)-l_0(x))y-(b(x)-b_0)y^2\right\}}-1\right|\left|\psi_k(y)\right|\left|g(x-y)\right|dy\right\}$$
$$+\left|\sum_{k_0(x)\leq k \leq k_1(x)}\left\{\int_{\TT}\psi_k(y)g(x-y)dy\right\}\right|=A(x)\:+\:B(x)\:.$$

Now for the first term, using \eqref{heur} and the small oscillation
of the exponential, we deduce that
$$A(x)\lesssim M_{\d}f(x)\:.$$

For the second term, as claimed initially, we remark that $B(x)$ is
the local version of \eqref{simplt}; to ``achieve" \eqref{concept}
we need to compare $B$ with some averages of the Hilbert transform.
Here, the main ingredient is
 \beq\label{v201}
\left|\sum_{k\leq
K}\psi_k(y)-2^{K-1}\int_{-2^{-K}}^{2^{-K}}R(y+z)dz\right|\lesssim\frac{2^{-K}}{\left|y\right|^2+2^{-2K}}\:,
\eeq where $R(y)=\sum_{k\in 10\N}\psi_{k}(y)$, $K\in\N$  and
$y\in\TT$.

Now, for $x\in I_j$ fixed, we conclude
$$ B(x)\lesssim\sup_{I\supset I_j}(\frac{1}{\left|I\right|}\int_{I}\left|(R*g)(y)\right|dy)+ \sup_{I\supset I_j}(\frac{1}{\left|I\right|}\int_{I}\left|g(y)\right|dy)$$
$$=M_{\d}(R*g)(x)+M_{\d}g(x)\:.$$

Finally, combining our estimates for $A$ and  $B$ and using the fact
that the operator $g\:\rightarrow\:R*g $ is bounded\footnote{This
comes from $\hat{R}\in L^{\infty}(\TT)$, which is an easy
consequence of the fact that the function $\psi$ is compactly
supported away from the origin and has mean zero.} on $L^2(\TT)$, we
conclude that
$$ \left\|T^{\p}f\right\|_2\lesssim \left\|M_{\d}\left(R*\left(M_{c_0}^{*}Q_{b_0}^{*}f\right)\right)\right\|_2
+\left\|M_{\d}f\right\|_2\lesssim \d^{1/2}\left\|f\right\|_2\:.$$

\end{proof}
\end{l1}

At this point we have learned how to estimate basic families of
tiles - having a simple geometric structure - for which we have
uniform control on the density factor. The next step (Lemmas 2 and
3) will be to understand the interaction between two such basic
families in the case in which we have no information about their
density factors, but we know that they are located in different
regions of the time-frequency plane.  (Here we will use the fact
that the geometric factor of pair $(P_1,P_2)$ is small whenever
$P_1$ and $P_2$ are not in the same family of tiles.)

Before presenting the lemmas, we will need several definitions.

\begin{d0}\label{sep}
Fix a number $\d\in(0,1)$. Let be $\p_1$ and $\p_2$ two trees with
(tops $\tilde{P}_1$ and $\tilde{P}_2$) representatives
$P_1=[\a_1,\o_1,I_1]$ and $P_2=[\a_2,\o_2,I_2]$ respectively; we say
that $\p_1$ and $\p_2$ are \emph{($\d$-)separated} if $\newline$
either $I_1\cap I_2=\es$ or else$\newline$i)
$P=[\a,\o,I]\in\p_1\:\:\&\:\:I\subseteq
I_2\:\:\:\:\Rightarrow\:\:\:\left\lceil
\Delta(P,P_2)\right\rceil<\d\:,$ $\newline$ii)
$P'=[\a',\o',I']\in\p_2\:\:\&\:\:I'\subseteq
I_1\:\:\:\:\Rightarrow\:\:\:\left\lceil
\Delta(P',P_1)\right\rceil<\d\:.$
\end{d0}
\noindent{\bf Notation:} Whenever we have two trees $\p_1$ and $\p_2$ as in
Definition \ref{sep} we will denote with
$x_{P_{1},P_{2}}^{i}:=x_{1,2}^{i}$ the abscissa of the
intersection point of $l_{P_1}$ with $l_{P_2}$. With this done, set
$\frac{w}{\min(|I_1|,|I_2|)}:=\frac{\left({\d}^{-1} \left\lceil
\Delta(P_1,P_2)\right\rceil\right)^{\frac{1}{2}}}{100}$ and define:
 \begin{itemize}
\item $I_s$ - the {\it separation interval (relative to the intersection)} of $\p_1$ and $\p_2$ by
$$I_s=[x_{1,2}^{i}-w,x_{1,2}^{i}+w]\cap\tilde{I}_1\cap\tilde{I}_2\:.$$

\item $I_c$ - the {\it ($\ep$-)critical
intersection interval} (between $\p_1$ and $\p_2$) by
$$I_c=3\d^{1/2-\ep}I_s$$ where $\ep$ is some small fixed positive real
number.
  \end{itemize}
\noindent{\bf Observation 5.} a) The two notions introduced above
can be regarded as indicators of how much the quadratic symmetry is
involved in the interaction of the two separated trees. Indeed, the
procedure of estimating terms like\footnote{Here $P_j\in\p_j$ with
$j\in\left\{1,2\right\}$.} \beq\label{ty}
 \left\langle {T_{P_1}}^{*}f , {T_{P_2}}^{*}g\right\rangle
\eeq will roughly obey the following scenario:

- if $\tilde{I}_{P_1}\cap I_s=\emptyset$ (or $\tilde{I}_{P_2}\cap I_s=\emptyset$)
then \eqref{ty} can be treated as in Fefferman's case, neglecting the
quadratic modulation

- else, guided by the results obtained in Section 5 (see Lemma 0),
we will split the integral in \eqref{ty} in two\footnote{We use the
fact that, for properly chosen $\ep$ and $\ep_0$, the
($\ep_0$-)critical intersection interval $I_{1,2}$ of the pair
$(P_1,P_2)$ is always included in $I_c$.}: the first (integrated
over the complement of $I_c$) will be treated as in the previous
case, while the second term (integrated over a set included in
$I_c$) will be placed into a collection of objects representing the
critical contribution of the quadratic symmetry.

b) Further, we will make use of two essential properties of our
above-defined intervals:

1) $\forall\:\:P\in\p_1\cup \p_2$ such that $x_{1,2}^{i}\in
5\tilde{I}_P$ we have $|I_P|>|I_s|$.

2) $\forall\:\:P\in\p_1\cup \p_2$ we have (for $\ep$ properly chosen) $|\tilde{I}_P\cap I_c|<{\d}^{\frac{1}{4}}|I_P|$.

\begin{l1}\label{sept}
Let be $\left\{\p_j\right\}_{j\in\left\{1,2\right\}}$ two separated
trees with tops $P_j=\left[\a_j,\o_j,I_0\right]$,
$j\in\left\{1,2\right\}$. Then, for any $f,\:g\in L^{2}(\TT)$ and
$n\in \N$, we have that \beq\label{v21} \left|\left\langle
{T^{\p_1}}^*f,{T^{\p_2}}^*g\right\rangle\right|\lesssim_{n}{\d}^n\left\|f\right\|_{L^{2}(\tilde{I}_0)}\left\|g\right\|_{L^{2}(\tilde{I}_0)}+\left\|\chi_{I_c}{T^{\p_1}}^*f\right\|_2\left\|\chi_{I_c}{T^{\p_2}}^*g\right\|_2\:.
\eeq (Remark: The first term in the right hand side of \eqref{v21}
expresses the result of the interaction\footnote{As if we were in
the linear phase case - see \cite{2}.} far from the intersection
point of the trees while the second one reflects the correction
needed for handling the quadratic case in the critical region $I_c$
.)
\begin{proof}

We start the proof of our lemma by making a partial (Whitney) dyadic
decomposition of the real axis with respect to the point
$x_{1,2}^{i}$; more exactly (we may assume that $|I_s|=2^{-r}$ for
some $r\in \N$), let

\begin{itemize}
    \item $A_0=\left(-\infty,\:x_{1,2}^{i}-\frac{|I_0|}{2}\right)\cup \left(x_{1,2}^{i}+\frac{|I_0|}{2},\:\infty\right)\\$

    \item $A_1=\left[x_{1,2}^{i}-\frac{|I_0|}{2},\:x_{1,2}^{i}-\frac{|I_0|}{4}\right)\cup \left(x_{1,2}^{i}+\frac{|I_0|}{4},\:x_{1,2}^{i}+\frac{|I_0|}{2}\right]\\$

 \item $A_2=\left[x_{1,2}^{i}-\frac{|I_0|}{4},\:x_{1,2}^{i}-\frac{|I_0|}{8}\right)\cup \left(x_{1,2}^{i}+\frac{|I_0|}{8},\:x_{1,2}^{i}+\frac{|I_0|}{4}\right] \\$

 \item $A_k=\left[x_{1,2}^{i}-|I_s|,\:x_{1,2}^{i}+|I_s|\right]\:.$
\end{itemize}

For $j\in\left\{1,2\right\}$ define the following sets:

$\newline S_{j,k}=\left\{P\in\p_j\:|\:x_{1,2}^{i}\in
5\tilde{I}_{P}\right\}\:\:\:\:\:\:\:\:\:\:\:\p'_j:=\p_j-S_{j,k}$
$\newline S_{j,0}=\left\{P\in\p'_j\:|\:{I_P}^{*}\cap
A_0\not={\es}\:\:\:\&\:\:\:|\tilde{I}_P|\leq\frac{|I_0|}{4}\right\}$
$\newline S_{j,l}=\left\{P\in\p'_j\:|\:{I_P}^{*}\cap
A_l\not={\es}\:,\:P\notin
S_{j,l-1}\:\:\:\&\:\:\:|\tilde{I}_P|\leq\frac{|A_l|}{3}\right\}\:\:\forall\:\:l\in\left\{1,..k-1\right\}\:.$

 With these notations it is clear that $\left\{
S_{j,l}\right\}_{l=1}^k$ form a partition of $\p_j $.

 Now setting
$T_{j,l}^{*}=\sum_{P\in S_{j,l}}{T_P}^{*}$ we obtain \beq\label{v22}
\left\langle{T^{\p_1}}^*,\:{T^{\p_2}}^*
\right\rangle=\sum_{n,l=0}^{k}\left\langle{T_{1,l}}^*,\:{T_{2,n}}^*
\right\rangle\:. \eeq

Now let $A_{-1}=\es$, $A_{k+1}=\es$ and
$\tilde{A}_{l}=^{def}\bigcup_{{P\in S_{j,l}}\atop{j\in
\left\{1,2\right\}}}\tilde{I}_P$; then
$supp\:{T_{j,l}}^{*}\subset\tilde{A}_{l}$ with
$\tilde{A}_{l}\subseteq A_{l-1}\cup A_{l}\cup
A_{l+1}\:\:\:\forall\:\:\:l\in\left\{0,..k-2\right\}\:\:\&\:\:\:\tilde{A}_{k-1}\cap\frac{1}{3}A_k={\es}$.
Consequently, to estimate $ \eqref{v22}$, we need to study the
following expressions\footnote{The terms of the form $\left\langle
{T_{1,l}}^*f,{T_{2,l-1}}^*g\right\rangle$ or
$\left\langle{T_{1,l}}^*f,{T_{2,l+1}}^*g\right\rangle$ have a
similar treatment.} \beq\label{T12}
\begin{array}{rl}

&U:=\sum_{l=0}^{k-1}\left\langle{T_{1,l}}^*f,{T_{2,l}}^*g
\right\rangle\:\:\:\:\:\:
Y:=\sum_{l=0}^{k-1}\left\langle{T_{1,k}}^*f,{T_{2,l}}^*g \right\rangle\:\\
&Z:=\sum_{l=0}^{k-1}\left\langle{T_{1,l}}^*f,{T_{2,k}}^*g
\right\rangle\:\:\:\:\:\:
V:=\left\langle{T_{1,k}}^{*}f,T_{2,k}^{*}g\right\rangle
\end{array}
\eeq

We concentrate now on the first term $U$.

For the beginning we will introduce several useful tools. For
$j\in\left\{1,2\right\}$, let $l_{P_{j}}(x)=l_j(x)=c_j+2xb_j$  and
$\newline d_{j,l}=\min\left\{|I|| P=[\a,\o,I]\in\p_j\:\&\:P\in
S_{j,l-1}\cup S_{j,l}\cup S_{j,l+1}\right\}$ (here $S_{j,-1},$
 $\\\:S_{j,k+1}:=\es$); also, define a real-valued function $\f\in C_{0}^{\infty}(\R)$ with the following properties:
 \beq\label{fprop}
 \eeq
 \begin{itemize}
 \item $ supp\:\f\subset\left\{\frac{1}{4}\leq |x|\leq \frac{1}{2}\right\}$
\item $ \f\:is\: even $
\item $ |\hat{\f}(\xi)-1|\lesssim_{n} |\xi|^n \:\:\:\:\forall\:|\xi|\leq 1\:\:and\:n\: big\:enough $
\item $ |\hat{\f}(\xi)|\lesssim_{n} |\xi|^{-n}\:\:\:\:\forall\:|\xi|\geq 1$
    \end{itemize}

Now, for $j\in\left\{1,2\right\}$ and $l\in\left\{1,..k-1\right\}$
set
$$\f_{j,l}(x)=(\d^{1/3}d_{j,l})^{-1}\f((\d^{1/3}d_{j,l})^{-1}x)$$
and define the operators
$$ \tilde{\f}_{j,l}\::L^2(\R)\longrightarrow L^2(\R)\:\:\:\: by\:\:\tilde{\f}_{j,l}f=\f_{j,l}*f $$ and
$$\Phi_{j,l}\::L^2(\R)\longrightarrow L^2(\R)\:\:\:\:by\:\:\Phi_{j,l}=M_{c_j}Q_{b_j}\tilde{\f}_{j,l}{Q_{b_j}}^{*}{M_{c_j}}^{*}\:.$$
Remark that $\tilde{\f}_{j,l}\:\:\textrm{and}\:\:\Phi_{j,l}$ are
self-adjoint for all $j\in\left\{1,2\right\}$ and
$l\in\left\{1,..k-1\right\}$.

Our first aim is to prove the following:
\begin{claim}
For$\: j\in\left\{1,2\right\}$ , $l\in\left\{1,..k-1\right\}$ and
$n\in \N$, decomposing ${T_{j,l}}^{*}$ as \beq\label{dec}
{T_{j,l}}^{*}f=\Phi_{j,l}{T_{j,l}}^{*}f\:+\:\Omega_{j,l}f \eeq we
have \beq \label{om} \left\|\Omega_{j,l}\right\|_2\lesssim_{n}\d^{n}
\eeq and \beq \label{fi} \left|\left\langle\Phi_{1,l}{T_{1,l}}^{*}f
,\Phi_{2,l}{T_{2,l}}^{*}g\right\rangle\right|\lesssim_{n}\d^{n}\left\|f\right\|_2\left\|g\right\|_2\:.
\eeq
\end{claim}

$\newline${\it Proof of the CLAIM} $\newline$

Suppose that $d_{2,l}\leq d_{1,l}$; now, since $\Phi_{j,l}$ are
self-adjoint, for showing \eqref{fi} it is enough to prove \beq
\label{om1}
\left|\left\langle\Phi_{2,l}\Phi_{1,l}\overbrace{{T_{1,l}}^{*}f}^{v}
,u\chi_{\tilde{A}_l}\right\rangle\right|\lesssim_{n}\d^{n}\left\|v\right\|_2\left\|u\right\|_2\:.
\eeq

Fix now $x\in \tilde{A}_l$; then
$$|\Phi_{2,l}\Phi_{1,l}v(x)|=|\tilde{\f}_{2,l}{Q_{b_2}}^{*}{M_{c_2}}^{*}M_{c_1}Q_{b_1}\tilde{\f}_{1,l}{Q_{b_1}}^{*}{M_{c_1}}^{*}v(x)|\leq $$ $$\int\left|v(s)\right|\left|\underbrace{\int{\f}_{2,l}(x-y){\f}_{1,l}(y-s)e^{i[(b_1-b_2)y^{2}+(c_2-c_1)y]}\:dy}_{\K_l(x,s)}\right|ds\:.$$

 Now making the change of variable $y=x-t\d^{1/3}d_{2,l}$ we deduce
$$|\K_l(x,s)|=$$
$$\left|\int_{\R}e^{i\left\{\d^{1/3}d_{2,l}t(l_2(x)-l_1(x))-(b_2-b_1)t^{2}(\d^{1/3}d_{2,l})^{2}\right\}}\f(t)\f_{1,l}(x-s-\d^{1/3}d_{2,l}t)dt\right|\:.$$

 Consequently, we need to estimate an expression of the form
$$I_{\a,\b}(\phi)=\left|\int_{\R}e^{i(\a t+\b t^{2})}\phi(t)dt\right|$$ where
$ \phi\in C_{0}^{\infty}(\R)\:$ with $\:\operatorname{supp}\:
\phi\subseteq\left\{\frac{1}{4}\leq|t|\leq\frac{1}{2}\right\}$, $
\a=\d^{1/3}d_{2,l}(l_2(x)-l_1(x))\:\:\textrm{and}\:\:\b=-(\d^{1/3}d_{2,l})^2
(b_2-b_1)\:.$

Now, since the trees are separated, and since $x\in\tilde{A}_l$, we
have that
$$\inf_{x\in\tilde{A}_l}\left|l_2(x)-l_1(x)\right|\geq \d^{-1}{d_{2,l}}^{-1}\:\:\:\Rightarrow\:\:\left|\a\right|\geq\d^{-2/3}\:.$$

Since for
$r(t)=t+\frac{\b}{\a}t^{2}\:\:\:\Rightarrow\:\:\:\left|r'(t)\right|\geq
1-\frac{|2\b|}{|\a|}t\geq 1-\frac{\d^{1/3}d_{2,l}}{2|A_l|}>0$
$\newline$we can apply the (non-)stationary phase method and deduce that
$$\left|I_{\a,\b}(\phi)\right|\lesssim_{n}{\a}^{-n}\:\:\:\forall\:\:n\in \N\:.$$

As a consequence, we have that
$$|\K_l(x,s)|\lesssim\d^n\underbrace{ (\d^{1/3}d_{1,l})^{-1}\chi_{\left\{|t|\leq2\d^{1/3}d_{1,l}\right\}}(x-s)}_{u_{1,l}(x-s)}\:\:\:\:\:\:\:\:\:\:\forall\:\:x\in\tilde{A}_l$$
$$\Rightarrow\:\:\:\:\left|\Phi_{2,l}\Phi_{1,l}\:v(x)\right|\lesssim\d^{n}(u_{1,l}*|v|)(x)\lesssim\d^{n}Mv(x)\:,$$ so \eqref{om1} holds.

We now discuss the expression
$$\Omega_{j,l}f ={T_{j,l}}^{*}f\:-\:\Phi_{j,l}{T_{j,l}}^{*}f\:. $$

Keeping in mind the fact that $Q_{b_j},M_{c_j}$ are unitary we have
the following chain of equalities:
$$\left\|\Omega_{j,l}f\right\|_2=\left\|{T_{j,l}}^{*}-\Phi_{j,l}{T_{j,l}}^{*}f\right\|_2=$$
$$\left\|{Q_{b_j}}^{*}{M_{c_j}}^{*}{T_{j,l}}^{*}f-\tilde{\f}_{j,l}{Q_{b_j}}^{*}{M_{c_j}}^{*}{T_{j,l}}^{*}f\right\|_2=^{f:=M_{c_j}Q_{b_j}h}$$
$$\left\|{Q_{b_j}}^{*}{M_{c_j}}^{*}{T_{j,l}}^{*}M_{c_j}Q_{b_j}h-\tilde{\f}_{j,l}{Q_{b_j}}^{*}{M_{c_j}}^{*}{T_{j,l}}^{*}M_{c_j}Q_{b_j}h\right\|_2\:.$$

Denote with
$\T_{j,l}^{*}={Q_{b_j}}^{*}{M_{c_j}}^{*}{T_{j,l}}^{*}M_{c_j}Q_{b_j}$;
then $\T_{j,l}={Q_{b_j}}^{*}{M_{c_j}}^{*}T_{j,l}M_{c_j}Q_{b_j}$ and
since
$\left\|\Omega_{j,l}\right\|_2=\left\|\Omega_{j,l}^{*}\right\|_2$ we
have that
$$\left\|\Omega_{j,l}\right\|_{2}=\left\|\T_{j,l}-\T_{j,l}\tilde{\f}_{j,l}\right\|_{2}\:.$$

Now fixing a tile $P=[\a,\o,I]\in S_{j,l}$ (with $|I|=2^{-k}$ for
some $k\in \N$) we set
$$\T_{j,l}^{P}={Q_{b_j}}^{*}{M_{c_j}}^{*}T_{P} M_{c_j}Q_{b_j}\:.$$

 Then, for $x\in E(P)\subseteq I$ we have
$$\left|\T_{j,l}^{P}h(x)-\T_{j,l}^{P}\tilde{\f}_{j,l}h(x)\right|=$$
$$\left|\int_{\TT}e^{i\left\{y(l_x(x)-l_j(x))-(b(x)-b_j)y^2\right\}}{\f}_{k}(y)[h-\f_{j,l}*h](x-y)dy\right|\leq$$
$$\int_{\TT}\left|(h\chi_{\tilde{A}_l})(x-y)\right|\left|r_{x}^{P}(y)-(r_{x}^{P}*{\f}_{j,l})(y)\right|dy$$
where  \beq \label{r}
 r_{x}^{P}(y)=e^{i\left\{y(l_x(x)-l_j(x))-(b(x)-b_j)y^2\right\}}{\f}_{k}(y)\:.
 \eeq

Our next step is to provide an $L^\infty$ bound on the expression
$$r_{x}^{P}(y)-(r_{x}^{P}*\f_{j,l})(y)=\int_{r}\hat{r_{x}^{P}}(\xi)(1-\hat{\f_{j,l}}(\xi))e^{i\xi
y}dy\:.$$

 For this we write
$$\left|\widehat{r_{x}^{P}}(\xi)\right|=\left|\int_{r}e^{i\left\{\frac{s}{2^k}[(l_x(x)-l_j(x))-\xi]-\frac{s^2}{2^{2k}}(b(x)-b_j)\right\}}\f(s)ds\right|$$
 and observe that $|l_x(x)-l_j(x)|\lesssim|\o_{P}|\lesssim 2^{k}$; from this, since $|\xi|\gtrsim 2^k$, we can apply the method of (non-)stationary phase to obtain
 \beq \label{rhat}
\left|\widehat{r_{x}^{P}}(\xi)\right|\lesssim
\left(1+\frac{|\xi|}{2^k}\right)^{-n}\:\:\:\:\:\:\forall\:\:n\in
\N\:\:\:\&\:\:\:x\in E(P)\:. \eeq

Using now  \eqref{rhat} together with \eqref{fprop}, we deduce that
$$\int_{r}\left|\widehat{r_{x}^{P}}(\xi)(1-\widehat{\f_{j,l}}(\xi))\right|d\xi\lesssim\int_{|\xi|\leq(\d^{1/3}d_j)^{-1}}\left(1+\frac{|\xi|}{2^k}\right)^{-n-2}({\d}^{1/3}d_{j,l})^n |\xi|^{n} d\xi\:+$$
$$\int_{|\xi|>(\d^{1/3}d_j)^{-1}}\left(1+\frac{|\xi|}{2^k}\right)^{-n-2} d\xi\lesssim 2^k(2^{k}\d^{1/3}d_{j,l})^n \:.$$

 As a consequence we have that for any $P\in\p_j$ and $x\in E(P)\subseteq I_P$
$$\left|r_{x}^{P}(y)-(r_{x}^{P}*\f_{j,l})(y)\right|\lesssim (\d^{1/3}2^{k}d_{j,l})^{n}2^{k}\chi_{[-2^{-k},2^{-k}]}(y)\:.$$

Then denoting
$\r_{j,l}(y)=\sum_{2^k\leq(d_{j,l})^{-1}}(\d^{1/3}2^{k}d_{j,l})^{n}2^{k}\chi_{[-2^{-k},2^{-k}]}(y)$
we obtain $\left\|\r_{j,l}\right\|_1\lesssim (\d)^{n/3}$, and so
$$\left\|\T_{j,l}h-\T_{j,l}\tilde{\f}_{j,l}h\right\|_{2}^{2}\lesssim\int_{\TT}\left(\int_{\TT}|h\chi_{\tilde{A}_l}(x-y)|\r_{j,l}(y)dy\right)^{2}dx\leq^{C-S}$$
$$\left\|\r_{j,l}\right\|_{1}^{2}\int_{\TT}|h\chi_{\tilde{A}_l}|^{2}\lesssim \d^{2n}\int_{\tilde{A}_l}|h|^{2}\:.$$

Consequently, we have shown that
$$\left\|\Omega_{j,l}^{*}h\right\|_{2}\lesssim\d^{n}\left(\int_{\tilde{A}_l}|h|^{2}\right)^{1/2}\:,$$
ending the proof of our claim.$\newline$

Now, reformulating the previous statements, we have
$\:\forall\:\:\:l\in\left\{0,..k-1\right\}$
$$\left\{\:{\left\|\Omega_{j,l}^{*}(f)\right\|_2\:,\:\left\|\Omega_{j,l}(f)\right\|_2\lesssim_{n}\d^{n}}\left(\int_{\tilde{A}_l}|f|^{2}\right)^{1/2}\atop{\left|\left\langle\Phi_{1,l}{T_{1,l}}^{*}f ,\Phi_{2,l}{T_{2,l}}^{*}g\right\rangle\right|\lesssim_{n}\d^{n}\left(\int_{\tilde{A}_l}|f|^{2}\right)^{1/2}\left(\int_{\tilde{A}_l}|g|^{2}\right)^{1/2}}\right.\:.$$
(Remark that
$\operatorname{supp}\:T_{j,l}\:,\operatorname{supp}\:T_{j,l}^{*},\operatorname{supp}\:\Omega_{j,l}^{*},\:\operatorname{supp}\:\Omega_{j,l}\subseteq\tilde{A}_l$\:.)

Now since        $$\left\langle{T_{1,l}}^{*}f
,{T_{2,l}}^{*}g\right\rangle=\left\langle\Phi_{1,l}{T_{1,l}}^{*}f
,\Phi_{2,l}{T_{2,l}}^{*}g\right\rangle+\left\langle\Omega_{1,l}f
,\Phi_{2,l}{T_{2,l}}^{*}g\right\rangle+$$
$$\left\langle\Phi_{1,l}{T_{1,l}}^{*}f ,\Omega_{2,l}g\right\rangle+\left\langle\Omega_{1,l}f ,\Omega_{2,l}g\right\rangle\:,\\$$ we have that
$$|U|=\left|\sum_{l=0}^{k-1}\left\langle{T_{1,l}}^*f,{T_{2,l}}^*g \right\rangle\right|\leq\sum_{l=0}^{k-1}\left|\left\langle\Phi_{1,l}{T_{1,l}}^{*}f ,\Phi_{2,l}{T_{2,l}}^{*}g\right\rangle\right|+$$
$$\sum_{l=0}^{k-1}\left|\left\langle\Omega_{1,l}f ,\Phi_{2,l}{T_{2,l}}^{*}g\right\rangle\right|+\sum_{l=0}^{k-1}\left|\left\langle\Phi_{1,l}{T_{1,l}}^{*}f ,\Omega_{2,l}g\right\rangle\right|+\sum_{l=0}^{k-1}\left|\left\langle\Omega_{1,l}f ,\Omega_{2,l}g\right\rangle\right|\lesssim$$
$$\d^{n}\sum_{l=0}^{k-1}\left(\int_{\tilde{A}_l}|f|^{2}\right)^{1/2}\left(\int_{\tilde{A}_l}|g|^{2}\right)^{1/2}+\sum_{l=0}^{k-1}\left\|\Omega_{1,l}(f)\right\|_2\left(\int_{\tilde{A}_l}|g|^{2}\right)^{1/2}+$$
$$\sum_{l=0}^{k-1}\left\|\Omega_{2,l}(g)\right\|_2\left(\int_{\tilde{A}_l}|f|^{2}\right)^{1/2}+\sum_{l=0}^{k-1}\left\|\Omega_{1,l}(f)\right\|_2\left\|\Omega_{2,l}(g)\right\|_2\lesssim^{C-S}\d^{n}\left\|f\right\|_2\left\|g\right\|_2\:.$$

The terms $Y$ and $Z$ can be treated similarly; we leave these
details for the reader.

 Now, it remains to estimate the term
$$V=\int_{\TT}{T^{*}_{1,k}}f\:{\overline{T^{*}_{2,k}g}}\:.$$

Setting $B_c=\TT\setminus\frac{1}{3}I_c$ we have that
$$V=
\left\langle
\chi_{1/3I_c}{T_{1,k}}^{*}f,\chi_{1/3I_c}{T_{2,k}}^{*}g\right\rangle+\left\langle
\chi_{B_c}{T_{1,k}}^{*}f,\chi_{B_c}{T_{2,k}}^{*}g\right\rangle$$
$$=A\:+\:B\:.$$

Clearly, only the second term requires some work; for this, we need
first to introduce some adapted tools: for
$\:j\in\left\{1,2\right\}$  and $\f$ as above, we define
$$\f_j(x)=\f(x)=(\d^{1/2}|I_s|)^{-1}\f((\d^{1/2}|I_s|)^{-1}x)\:\:\&\:\:\tilde{\f}_j(f)=\tilde{\f}(f)=\f*f$$
$$\Phi_{j}\::L^2(\R)\longrightarrow L^2(\R)\:\:\:\:\textrm{with}\:\:\Phi_{j}=M_{c_j}Q_{b_j}\tilde{\f_{j}}{Q_{b_j}}^{*}{M_{c_j}}^{*}\:,$$
and finally $\:\:\:\:\:\Omega_{j}f
:=\chi_{B_c}{T_{j,k}}^{*}f\:-\:\Phi_{j}\chi_{B_c}{T_{j,k}}^{*}f$.

 Then for $j\in\left\{1,2\right\}$, we have
$$\chi_{B_c}T_{j,k}^{*}f=\Phi_{j}\chi_{B_c}T_{j,k}^{*}f\:+\:\chi_{\frac{2}{3}I_c}\Omega_{j}f\:+\:\chi_{\:c_{\left(\frac{2}{3}I_c\right)}}\Omega_{j}f\:. $$

Using the facts: $\newline$i) for $x\in B_c$ we have
$|l_1(x)-l_2(x)|\gtrsim {\d}^{-\ep}|\operatorname{supp}\:\f|^{-1}$
$\newline$ii)
$\chi_{\:c_{\left(\frac{2}{3}I_c\right)}}\Omega_{j}f\:=\chi_{\:c_{\left(\frac{2}{3}I_c\right)}}{T_{j,k}}^{*}f\:-\:\chi_{c_{\left(\frac{2}{3}I_c\right)}}\Phi_{j}{T_{j,k}}^{*}f$
$\newline$we can repeat the previous arguments and obtain

$$\left\|\Phi_{j}*g\right\|_{2}\lesssim\left\|g\right\|_2$$
$$\left\|\chi_{\:{\frac{2}{3}I_c}}\Phi_{j}\chi_{B_c}{T_{j,k}}^{*}f\right\|_2\lesssim\left\|\chi_{I_c}{T_{j,k}}^*f\right\|_2=\left\|\chi_{I_c}{T^{\p_j}}^*f\right\|_2$$
$$\left\|\chi_{\:c_{\left(\frac{2}{3}I_c\right)}}\Omega_{j}f\right\|_2\lesssim\d^{n}\left\|f\right\|_2$$
$$\left|\left\langle\Phi_{1}\chi_{B_c}{T_{1,k}}^{*}f ,\Phi_{2}\chi_{B_c}{T_{2,k}}^{*}g\right\rangle\right|\lesssim_{n}\d^{n}\left(\int_{\tilde{I_0}}|f|^{2}\right)^{1/2}\left(\int_{\tilde{I_0}}|g|^{2}\right)^{1/2}\:.$$

Putting these relations together we conclude
$$\left|B\right|\lesssim_{n}{\d}^n\left(\int_{\tilde{I_0}}|f|^{2}\right)^{1/2}\left(\int_{\tilde{I_0}}|g|^{2}\right)^{1/2}+\left\|\chi_{I_c}{T^{\p_1}}^*f\right\|_2\left\|\chi_{I_c}{T^{\p_2}}^*g\right\|_2\:.$$

Since we trivially have
$$\left|A\right|\leq\left\|\chi_{I_c}{T^{\p_1}}^*f\right\|_2\left\|\chi_{I_c}{T^{\p_2}}^*g\right\|_2\:,$$
our proof is now complete.
\end{proof}
\end{l1}

\begin{d0}\label{nor}
A tree $\p$ with top-representative $P_0=[\a_0,\o_0,I_0]$ is called
\emph{normal} if $\newline$
$P=[\a,\o,I]\in\p\:\:\:\:\:\:\Rightarrow\:\:\:\:\:\:|I|\leq\frac{\d^{100}}{K}|I_0|\:\:\:\&\:\:\:\:dist(I,\partial
I_0)>20\frac{\d^{100}}{K}|I_0|$.
\newline(Here $K>10$ is some fixed constant and $ \partial I_0$ designates the
boundary of $I_0$.)
\end{d0}
\noindent{\bf Observation 6.} Notice that if $\p$ is a normal tree
then $supp\:{{T^{\p}}^{*}f}\subseteq \left\{x\in
I_0\:|\:dist(x,\partial I_0)>10\frac{\d^{100}}{K}|I_0|\right\}$.

\begin{d0}\label{row}
A \emph{row} is a collection $\p=\cup_{j\in \N}\p^{j}$ of normal
trees $\p^{j}$ with top-representatives
$P^{j}_{0}=[\a^j_{0},\o^j_0,I^j_{0}]$ such that the
$\left\{I^j_0\right\}$ are pairwise disjoint.
\end{d0}

\begin{l1}\label{rt}
Let $\p$ be a row as above, let $\p'$ be a tree with
top-representative $[\a_{0}',\o_{0}',I_{0}']$ and suppose that
$\forall\: j\in \N$,$\:\:I_0^{j}\subseteq I_0'$ and $\p^{j},\:\p'$
are separated trees; for each $j$, denote by $I_c^{j}$ the critical
intersection interval between $\p^{j}$ and $\p'$.

Then for any $f,\:g\in L^{2}(\TT)$ and $n\in \N$ we have that
$$\left|\left\langle {T^{\p'}}^*f,{T^{\p}}^*g\right\rangle\right|\lesssim_{n}{\d}^n\left\|f\right\|_{2}\left\|g\right\|_{2}+\left\|\sum_{j}\chi_{I_c^{j}}{T^{\p'}}^*f\right\|_2\left\|\sum_{j}\chi_{I_c^{j}}{T^{\p^{j}}}^*g\right\|_2\:.$$
\begin{proof}

First, observe that it is enough to show that, for a fixed $j$, we
have
$$\left|\left\langle {T^{\p'}}^*f,{T^{\p^{j}}}^*g\right\rangle\right|\lesssim_{n}\d^n\left(\left\|M(Mf)\right\|_{L_2(I_0^j)}+\left\|M(M({T^{\p'}}^*f))\right\|_{L_2(I_0^j)}\right)\left\|g\right\|_{L_2(I_0^j)}$$
$$+\left\|\chi_{I_c^{j}}{T^{\p'}}^*f\right\|_2\left\|\chi_{I_c^{j}}{T^{\p^{j}}}^*g\right\|_2\:.$$

For simplicity, in what follows we will drop the index $j$.
Repeating now the procedures from the previous lemma, we define the
following objects: $
\left\{A_l\right\}_{\left\{l\in\left\{0,..k\right\}\right\}}$ - the
dyadic decomposition with respect to the (abscissa of the)
intersection point - $x^{i}$, $\p=\cup_{l=0}^{k}S_{l}$ - the
partition of the tree $\p$ in the well-localized (with respect to
the separation interval $I_s$) sets of tiles, and
$\left\{T_{l}^*\right\}_{\left\{l\in\left\{0,..k\right\}\right\}}$
the corresponding decomposition of ${T^{\p}}^{*}$ (so we have
$\\{T^{\p}}^{*}=\sum_{l=0}^{k}T_{l}^*$). Also, for
$$\\d_{l}=\min\left\{|I||P=[\a,\o,I]\in\p\cup\p'\:\:\&\:\:I\subseteq
A_{l-1}\cup A_l\cup A_{l+1},\:I\subseteq I_0\right\}$$
 define $\f_{l}\:,\:\tilde{\f_{l}}\:,\:\Phi_{l}$ and $\Omega_{l}$ as before. Finally, set
$$\Phi_{l}'=M_{c'}Q_{b'}\tilde{\f_{l}}{Q_{b'}}^{*}{M_{c'}}^{*}\:\:\&\:\:\Omega_{l}'f={T^{\p'}}^*f-\Phi_{l}'\:{T^{\p'}}^*f\:.$$

Then, for $l\in\left\{0,..k-1\right\}$, we have \beq \label{rowdec}
\begin{array}{rl}
\left\langle {T^{\p'}}^*f,{T_{l}}^*g\right\rangle=\left\langle
\Phi_{l}'{T^{\p'}}^*f,\Phi_{l}{T_{l}}^*g\right\rangle+\\
\left\langle\Phi_{l}'{T^{\p'}}^*f,\Omega_{l}g
\right\rangle+\left\langle\Omega_{l}'f,{T_{l}}^*g \right\rangle\:.
\end{array}
\eeq

Now using the following relations (see the previous lemma):
$$\left|\Phi_{l}'\Phi_{l}h(x)\right|\lesssim_{n}\d^n \left(u_{l}*|h|\right)(x)$$
$$\left|\Omega_{l}^{*}h(x)\right|\lesssim\left(\r_{l}*|h|\right)(x)\:\:$$
 with $\left\|\r_{l}\right\|_1\lesssim (\d)^{n/3}$, where
 $$\:\r_{l}(y)=\sum_{2^k\leq(d_{l})^{-1}}(\d^{1/3}2^{k}d_{l})^{n}2^{k}\chi_{[-2^{-k},2^{-k}]}(y)$$
and
$$u_{l}(x)=(d^{1/3}d_{l})^{-1}\chi_{\left\{|t|\leq2\d^{1/3}d_{l}\right\}}(x)\:,$$
$\newline$ we deduce (recall that $\p$ is a normal tree) for the
first two terms \beq \label{term1}
\begin{array}{rl}
\left|\left\langle \Phi_{l}'{T^{\p'}}^*f,\Phi_{l}{T_{l}}^*g\right\rangle\right|\lesssim_{n}\d^n\left\langle M\left({T^{\p'}}^*f\right),|T_{l}^*g|\right\rangle\lesssim\\
\d^n\left\|M({T^{\p'}}^*f)\right\|_{L^2(\tilde{A}_{l}\cap
I_0)}\left\|g\right\|_{L^2(\tilde{A}_{l}\cap I_0)}
\end{array}
\eeq and respectively, \beq \label{term2}
\begin{array}{rl}
\left|\left\langle\Phi_{l}'{T^{\p'}}^*f,\Omega_{l}g \right\rangle\right|\lesssim\left\langle \r_{l}*\left\{\chi_{\tilde{A}_{l}\cap I_0}M({T^{\p'}}^*f)\right\},|g|\right\rangle\lesssim\\
\d^n\left\|M({T^{\p'}}^*f)\right\|_{L^2(\tilde{A}_{l}\cap
I_0)}\left\|Mg\right\|_{L^2(\tilde{A}_{l}\cap I_0)}\:.
\end{array}
\eeq

 We now treat the last term of the right-hand side of
\eqref{rowdec}. For this, set first
$\p'_{l}:=\left\{P=[\a,\o,I]\in\p'|\:|I|\geq d_{l}\right\}$; then,
for $x\in supp\:T_l^{*}$, we have
$$                  \Omega_{l}'f(x)={T^{\p'}}^*f(x)-\Phi_{l}'{T^{\p'}}^*f(x)={T^{\p'_{l}}}^*f(x)-\Phi_{l}'{T^{\p'_{l}}}^*f(x)$$ and consequently we deduce

$$x\in supp\:T_l^{*}\:\:\Rightarrow\:\:\:\left|{\Omega_{l}'}^{*}h(x)\right|\lesssim\left(\r_{l}*|h|\right)(x)\:,$$
 so
\beq \label{term3}
\begin{array}{rl}

\left|\left\langle\Omega_{l}'f,{T_{l}}^*g
\right\rangle\right|\lesssim\left\langle |f|,\r_{l}*|{T_{l}}^*g|
\right\rangle\lesssim\d^n\left\|Mf\right\|_{L^2(\tilde{A}_{l}\cap
I_0)}\left\|g\right\|_{L^2(\tilde{A}_{l}\cap I_0)}\:.
\end{array}
\eeq

Now, adding the relations \eqref{term1} - \eqref{term3}, we obtain
$$\left|\left\langle {T^{\p'}}^*f,{T_{l}}^*g\right\rangle\right|\lesssim\d^n\left(\left\|Mf\right\|_{L^2(\tilde{A}_{l}\cap I_0)}+\left\|M({T^{\p'}}^*f)\right\|_{L^2(\tilde{A}_{l}\cap I_0)}\right)\left\|Mg\right\|_{L^2(\tilde{A}_{l}\cap I_0)}$$
and consequently from Cauchy-Schwarz we deduce \beq \label{rowe}
\begin{array}{rl}

&\left|\left\langle {T^{\p'}}^*f,{T^{\p}}^*g\right\rangle\right|\lesssim_{n}\\
&\d^n\left(\left\|Mf\right\|_{L^2(I_0)}+\left\|M({T^{\p'}}^*f)\right\|_{L^2(
I_0)}\right)\left\|Mg\right\|_{L^2(I_0)}+ \left|\left\langle
{T^{\p'}}^*f,T_{k}^{*}g\right\rangle\right|
\end{array}
\eeq where $T_{k}^{*}=\sum_{P\in S_k}T_{P}^{*}$ with
$S_k=\left\{P\in\p|x^{i}\in 5\tilde{I}_{P}\right\}\:.$

Now, for the last term of the right-hand side of \eqref{rowe}, we
argue as follows:

\noindent Case 1: $\:\:|I_s|\gtrsim\frac{\d^{100}}{K}|I_0|$

In this situation we have no tile $P=[\a,\o,I]\in\p$ such that
$100I\cap I_s\not=\es$, and consequently $T_{k}^*=0$, so we have
nothing to prove.

\noindent Case 2: $\:\:|I_s|\lesssim\frac{\d^{100}}{K}|I_0|$

Let be ${\p_{i}}'=\left\{P\in\p'|x^{i}\in 5\tilde{I}_{P}\right\}$,
$\grave{\p}=\p'\setminus\left(\p'_{n}\cup{\p_{i}}'\right)$ where
$\newline\p'_{n}=\left\{P=[\a,\o,I]\in\left\{\p'\setminus{\p_{i}}'\right\}|\:|I|\leq\frac{\d^{100}}{K}|I_0|\right\}$;
define ${T'_{i}}^{*}=\sum_{P\in {\p_{i}}'}T_{P}^{*}$; then obviously
$${T^{\p'}}^*f={T^{\p'_{n}}}^*f\:+\:{T^{\grave{\p}}}^*f\:+\:{T_{i}'}^{*}f\:.$$

For the first term, from Lemma 2, we deduce that
$$\left|\left\langle {T^{\p'_{n}}}^*f,T_k^{*}g\right\rangle\right|\lesssim_{n}{\d}^n\left\|f\right\|_{L^2(I_0)}\left\|g\right\|_{L^2(I_0)}\:.$$

Now using that
$P=[\a,\o,I]\in\grave{\p}\:\Rightarrow\:|I|\gtrsim|I_s|$ and
defining
$$\f(x)=\f(x)=(\d^{1/2}|I_s|)^{-1}\f((\d^{1/2}|I_s|)^{-1}x)\:\:\&\:\:\tilde{\f}(f)=\tilde{\f}(f)=\f*f$$
$$\Phi=M_{c}Q_{b}\tilde{\f}{Q_{b}}^{*}{M_{c}}^{*}\:\:\&\:\:\Phi'=M_{c'}Q_{b'}\tilde{\f}{Q_{b'}}^{*}{M_{c'}}^{*}$$
we can follow the general ideas presented above and show that
$$\left|\left\langle {T^{\grave{\p}}}^*f,T_{k}^{*}g\right\rangle\right|\lesssim_{n}{\d}^n\left(\left\|Mf\right\|_{L^2( I_0)}+\left\|M\left\{\M({T^{\p'}}^*f)\right\}\right\|_{L^2( I_0)}\right)\left\|Mg\right\|_{L^2(I_0)}$$
and
$$\left|\left\langle {T_{i}'}^{*}f,T_{k}^{*}g\right\rangle\right|\lesssim_{n}{\d}^n\left(\left\|Mf\right\|_{L^2( I_0)}+\left\|M\left\{\M({T^{\p'}}^*f)\right\}\right\|_{L^2( I_0)}\right)\left\|Mg\right\|_{L^2(I_0)}$$
$$+\left\|\chi_{I_c}{T^{\p'}}^*f\right\|_2\left\|\chi_{I_c}{T^{\p}}^*g\right\|_2\:,$$
where
$$\M({T^{\p'}}^*f)=^{def}\sup_{m\in N}\left|\sum_{{P=[\a,\o,I]\in\ \p'}\atop{|I|\geq 2^{-m}}}T_{P}^{*}f\right|\:.$$

So to summarize, we proved that
$$\left|\left\langle {T^{\p'}}^*f,{T^{\p}}^*g\right\rangle\right|\lesssim_{n}\d^n\left(\left\|Mf\right\|_{L^2( I_0)}+\left\|M\left\{\M({T^{\p'}}^*f)\right\}\right\|_{L^2( I_0)}\right)\left\|Mg\right\|_{L^2(I_0)}$$
$$+\left\|\chi_{I_c}{T^{\p'}}^*f\right\|_2\left\|\chi_{I_c}{T^{\p}}^*g\right\|_2\:.$$

Now the conclusion follows if we add the observation that \beq
\label{treemax} \M({T^{\p'}}^*f)\leq Mf\:+\:M({T^{\p'}}^*f)\:. \eeq

Indeed, we first see that \eqref{treemax} can be rewritten as \beq
\label{rot} \M({{\T}^{\p'}}^{*}f)\leq Mf\:+\:M({{\T}^{\p'}}^{*}f)
\eeq where, as usual,
${{\T}^{\p'}}^*:=Q_{b_0'}^{*}M_{c_0'}^{*}{T^{\p'}}^{*}M_{c_0'}^{*}Q_{b_0'}^{*}$
and $l'(x)=c_0'+2b_0'x$ is the central line associated with the top
of $\p'$.

Fix now $m\in \N$, $x\in[0,\:1]$  and define a function $\phi\in
C_{0}^{\infty}(\R)$ with the following properties:
\beq\label{fiprop} \eeq
 \begin{itemize}
\item $ supp\:\phi\subset\left\{|x|\leq 2\right\}$
\item $ \int_{\R}\phi=1\:\:\:\:\:\&\:\:\:\phi\geq 0\:.$
    \end{itemize}

Let be $J$ the dyadic interval having the properties $x\in J$ and
$|J|=2^{-m}$; set
$\phi_{J}(x):=|J|^{-1}\phi\left(|J|^{-1}(x-c(J))\right)$; we want to
estimate the expression
$$(A)=\left|\sum_{P\in\p'\atop{|I_P|\geq |J|}}{{\T}_{P}}^*f(x)-\int_{\R}\phi_{J}(s){{\T}^{\p'}}^*f(s)ds\right|\leq$$
$$\int_{\R}\phi_{J}(s)\left\{\sum_{P\in\p'\atop{|I_P|\geq |J|}}|{{\T}_{P}}^*f(x)-{{\T}_{P}}^*f(s)|\right\}ds \:+\:\left|\int_{\R}\phi_{J}(s)\left\{\sum_{P\in\p'\atop{|I_P|<|J|}}{{\T}_{P}}^*f(s)\right\}ds\right|$$
$$:=B+C\:.$$

We start by treating the first term; observe first, that with the
notation \eqref{r} we have (up to conjugation) that \beq\label{t}
{{\T}_{P}}^*f(x)=\int_{Pi}r_{y}^P(x-y)f(y)\chi_{E(P)}(y)dy\:. \eeq

Relying on this, we further have

$$\left|{{\T}_{P}}^*f(x)-{{\T}_{P}}^*f(s)\right|\leq \int_{\TT}|r_{y}^P(x-y)-r_{y}^P(s-y)||f(y)|\chi_{E(P)}(y)dy\:.$$

Using now relation \eqref{rhat} (for $|I_P|=2^{-k}$) we deduce
\beq\label{meanval}
|r_{y}^P(x-y)-r_{y}^P(s-y)|=\left|\int_{\R}\hat{r}_{y}^P(\eta)\left(e^{i\eta(x-y)}-e^{i\eta(s-y)}\right)d\eta\right|\lesssim
\eeq
$$\int_{\R}\left(1+\frac{\eta}{2^k}\right)^{-n} \eta|x-s|d\eta\lesssim 2^{2k}|J|\:,$$
where we used the fact that $x,\:s\in 5J$ and $y\in E(P)$.

From the previous relations, we conclude
$$B\lesssim\int_{\R}\phi_{J}(s)\left\{\sum_{P\in\p'\atop{J\subseteq 3\tilde{I}_P}}(|I_P|^{-1}|J|)\int_{\R}|f(y)|\frac{\chi_{I_P}(y)}{|I_P|}dy\right\}ds\lesssim Mf(x)\:.$$

For the second term $C$ we use the fact that $\psi$ has the mean
zero property. Indeed, we have
$$C=\left|\sum_{P\in\p'\atop{|I_P|<|J|}}\int_{\R}f(y)\chi_{E(P)}(y)\left(\int_{\R}\phi_{J}(s)r_{y}^P(s-y)ds\right)dy\right|\lesssim$$
$$\sum_{P\in\p'\atop{\tilde{I}_P\subset10J}}\int_{\R}|f(y)|\chi_{E(P)}(y)\left(\int_{\R}|\hat{\phi}_{J}(\eta)\hat{r}_{y}^P(\eta)|d\eta\right)dy\:.$$

Now, for a fixed $y$, we argue as follows:
$$\sum_{P\in\p'\atop{\tilde{I}_P\subset10J}}\chi_{E(P)}(y)\left(\int_{\R}|\hat{\phi}_{J}(\eta)\hat{r}_{y}^P(\eta)|d\eta\right)\lesssim$$
$$\chi_{10J}(y)\int_{\R}|\hat{\phi}_{J}(\eta)|\left(1+\sum_{k>m}2^{-k}|\eta|\right)d\eta \lesssim|J|\:\chi_{10J}(y)\:.$$
Consequently,
$$C\lesssim Mf(x)$$
and replacing the bounds for $B$ and $C$ in $(A)$, we conclude
$$\left|\sum_{P\in\p'\atop{|I_P|\geq |J|}}{{\T}_{P}}^*f(x)\right|\lesssim\left|\int_{\R}\phi_{J}(s){{\T}^{\p'}}^*f(s)ds\right|+Mf(x)$$
which is what we needed for \eqref{treemax} to hold.
\end{proof}
\end{l1}

Finally, we combine the previous results to prove that we can
control the $L^2$ norm of the operator associated to a forest.

$\newline${\bf Proof of Proposition 2}

$\newline$

Define $F=\bigcup_{j}\left\{x\in
I_j\:|\operatorname{dist}(x,\partial
I_j))\leq100\frac{\d^{100}}{K^2}|I_j|
\right\}=^{def}\bigcup_{j}F_j\:.$

We will estimate our operator only on the complement of this set.
This is safe since we can control the measure of the excised set as
follows:
$$|F|\leq\sum_{j}|F_j|\lesssim\sum_j |I_j|\frac{\d^{100}}{K^2}\lesssim\frac{\d^{50}}{K}\:,$$
where the last inequality it is derived from hypothesis 3).

 Now, on the $F^{c}$, we intend to use the previous
estimates obtained in Lemma 3, but before this, we are forced to
create enough space\footnote{Here it is essential that our trees are
``centered" - see Observation 3 a).} to apply the separation
results. Consequently, we start by removing few tiles\footnote{In
the following procedure, we will assume that there is no tree $\p_j$
having two tiles with same time interval; if this is not the case,
then we must have (for some j) the situation $P\in\p_j$ and $P_u\in
\p_j$ (or $P_l\in\p_j$), in which case we take the union of these
two tiles and consider it as a single tile - renamed $P$.} from each
tree $\p_j$.

Let be $\p=\cup_{j}\p_j$; for $M=\log{(K^{100}\d^{-100})}\:$ denote
$$\p^{+}:=\left\{P\in\p\:|\:\operatorname{there\:is\:no\:chain\:}
 P<P_1<...<P_{M}\:\operatorname{with\:all\:}P_j\in\p\right\}$$
and
$$\p^{-}:=\left\{P\in\p\:|\:\operatorname{there\:is\:no\:chain\:}
 P_1<P_2<...<P_{M}<P\:\operatorname{with\:all\:}P_j\in\p\right\}\:.$$

Now, it is easy to see that each such set can be split into at most
$M$ subsets with no two comparable tiles inside the same subset.
Consequently, using Proposition 1, we deduce that\footnote{As
mentioned in Section 3, $\eta$ may change from line to line.}
$$\left\|T^{\p^{+}}\right\|_2\:,\:\left\|T^{\p^{-}}\right\|_2\lesssim M{\d}^{\eta}\lesssim{\d}^{\eta}\log{K}\:.$$
We remove all the above mentioned sets from our collection $\p$ and
decompose this new set as follows:
$$\p=\cup_{j}\p^{0}_j\:\:\operatorname{where}\:\: \p^{0}_j=\p_j\cap\p\:.$$

Now this modified collection $\p$ behaves much better than the
initial one; indeed, we have
$\newline\:1)\:\forall\:P=[\a,\o,I]\in\p^{0}_j$,
$\:|I|\leq\frac{\d^{100}}{K^{100}}|I_j|$
$\newline\:2)\:\forall\:j\not=k$, the trees $\p^{0}_j$ and
$\p^{0}_k$ are $\d'$-separated where
$\d'=\frac{\d^{100}}{K^{100}}\:.$

 Moreover, if we split each
$\p^{0}_j=\p^{N}_j\cup\p^{C}_j$, with
$$\p^{C}_j=^{def}\left\{P=[\a,\o,I]\in\p^{0}_j\:|\:I\subseteq F_j\right\}\:,$$
we conclude that $\:\left\{\p^{N}_j\right\}_{j}$ represents a
collection of normal, $\d'$-separated trees, while for the remaining
parts of the trees we have the relation
$$\operatorname{supp}\:T^{\p^{C}_j}\subset F_j\:.$$

 Consequently, on $F^c$ we have that
$$T^{\p}f\:=\:\sum_{j}T^{\p^{N}_j}f\:,$$
and so our conclusion reduces to \beq\label{enough}
\left\|\sum_{j}T^{\p^{N}_j}f\right\|_2\lesssim
\d^{\eta}\log{K}\left\|f\right\|_2\:. \eeq

Now we are ready to apply the results from Lemma 3. We start by
dividing $\bigcup_{j}\p^{N}_j$ into a union of at most $K\d^{-2}$
rows, $\r_1\:,\:\r_2\:,...\r_{K\d^{-2}}$. This is done by using an
easy maximal argument: choose from $\left\{I_j\right\}_j $ a
collection of maximal (disjoint) dyadic intervals - call it $r_1$;
after that, erase the set $r_1$ from the previous collection and
repeat the same procedure with the remaining one obtaining a new set
$r_2$; due to condition 3) in our hypothesis, we know that this
procedure will end in at most $K\d^{-2}$ steps; now take $\r_j$ to
be the set containing all trees that have their top inside the set
$r_j$. Now, denoting by $T^{\r_j}$ the operator associated with
$\r_j$, we claim that $\\\newline
C1)\:\left\|T^{\r_{j}}\right\|_{\left\{2\rightarrow
2\right\}}\lesssim {\d}^{\frac{1}{2}}\:.$ $\newline
C2)\:\left\|{T^{\r_{k}}}^{*}T^{\r_{j}}\right\|_{\left\{2\rightarrow
2\right\}}=0\:\:\:\:\:$ if $\:\:\:k\not=j\:.$ $\newline
C3)\:\left\|T^{\r_{k}}{T^{\r_{j}}}^{*}\right\|_{\left\{2\rightarrow
2\right\}}\lesssim \left(\frac{\d}{K}\right)^{10}\:\:\:\:$  if
$\:\:\:k\not=j\:.$ $\\\newline$

If we accept this for the moment, then applying the Cotlar-Stein
Lemma we deduce that
$$\left\|\sum_{j}T^{\p^{N}_j}f\right\|_2=\left\|\sum_{j}T^{{\r}_j}f\right\|_2\lesssim \d^{\frac{1}{2}}\left\|f\right\|_2\:.$$

This last relation trivially implies \eqref{enough}, ending our proof.

 We now pass to the analysis of our claims; for C1), we
just remark that since $\r_j$ is a row all the trees that belong to
it are spatially disjoint, which together with Lemma 1 implies our
statement. C2) is trivial since for $k\not=j$, the operators
$T^{\r_{k}}$ and $T^{\r_{j}}$ live in disjoint parts of the unit
interval. The only interesting claim is C3). Fix $k_0,\:j_0$ and
suppose that $j_0<k_0$. To avoid working with double indices, we
will make the following notations: let $\left\{\A_j\right\}_j$ be
the tree-decomposition of $\r_{j_0}$ with top time intervals
$\left\{A_j\right\}_j$ and $\left\{\B_{k}\right\}_k$ the trees
corresponding to $\r_{k_0}$ and with top time intervals
$\left\{B_k\right\}_k$. Since $j_0<k_0$ (from the way we constructed
our rows) we have that $A_j\cap B_k\not=\es$ implies $B_k\subseteq
A_j$. Given this fact, we may assume that there exists
$\left\{n_l\right\}_l\subset \N$ a strictly increasing sequence of
natural numbers ($n_0=1$) such that
$$A_j\supseteq\bigcup_{n_{j-1}\leq k<n_j}B_k\:.$$

Now, from the fact that our trees are normal, we have
$$\left\langle {T^{\r_{j_0}}}^{*}f,{T^{\r_{k_0}}}^{*}g\right\rangle=\sum_{j\geq 1}\left\langle{T^{\A_j}}^{*}f,\sum_{n_{j-1}\leq k<n_j}{T^{\B_k}}^{*}g \right\rangle\:,$$
where here $f,\:g$ are two arbitrary functions in $L^2(\TT)$.

Now define $I^{j,k}_{c}$ to be the critical intersection interval
associated with the trees $\A_j$ and $\B_k$, and let
 $I^j_c=\bigcup_{n_{j-1}\leq k<n_j}I^{j,k}_{c}$. Then, applying Lemma 3 for a fixed $j$, we deduce
\beq \label{jk} \left|\left\langle{T^{\A_j}}^{*}f,\sum_{n_{j-1}\leq
k<n_j}{T^{\B_k}}^{*}g
\right\rangle\right|\lesssim_{n}(\d')^n\left\|f\right\|_{L^{2}(A_j)}\left\|g\right\|_{L^{2}(A_j)}+
\eeq
$$\left\|\chi_{I_c^{j}}{T^{\A_j}}^*f\right\|_2\left\|\sum_{n_{j-1}\leq k<n_j}\chi_{I_c^{j,k}}{T^{\B_k}}^*g\right\|_2\:,$$
where we remind that the trees $\A_j$ and $\B_k$ are $\d'$-separated
with $\d'=\frac{\d^{100}}{K^{100}}$.

Now, using the relations\footnote{See Observation 5 b).} \beq
\label{inter1} \: |I^{j}_{c}\cap \tilde{I}_{P}|\leq
\frac{{\d}^{25}}{K^{25}}|I_P|\:\:\:\:\:\:\:\:\:\forall\:P\in\A_{j}
\eeq \beq \label{inter2} \: |I^{j,k}_{c}\cap \tilde{I}_{P}|\leq
\frac{{\d}^{25}}{K^{25}}|I_P|\:\:\:\:\:\:\:\:\:\forall\:P\in\B_{k}
\eeq
 together with Lemma 4 below and the fact that all trees involved are normal, we have that
$$\left\|\chi_{I_c^{j}}{T^{\A_j}}^*f\right\|_2\lesssim\left(\frac{\d}{K}\right)^{10}\left\|f\right\|_{L^{2}(A_j)}$$
and
$$\left\|\sum_{n_{j-1}\leq k<n_j}\chi_{I_c^{j,k}}{T^{\B_k}}^*g\right\|_2\lesssim\left(\frac{\d}{K}\right)^{10}\left\|g\right\|_{L^{2}(A_j)}\:.$$

Now replacing the last two relations in \eqref{jk}, we conclude
$$\left|\left\langle{T^{\A_j}}^{*}f,\sum_{n_{j-1}\leq k<n_j}{T^{\B_k}}^{*}g  \right\rangle\right|\lesssim \left(\frac{\d}{K}\right)^{10}\left\|f\right\|_{L^{2}(A_j)}\left\|g\right\|_{L^{2}(A_j)}\:,$$
which together with an easy orthogonality argument gives us relation
C3), completing our proof.

\begin{flushright}
$\Box$
\end{flushright}

Heuristically, the next result is a ``dual statement" of Lemma 1.

\begin{l1}\label{ed}
Let $\p\:$ be a tree with top-representative
$P_{0}=[\a_{0},\o_{0},I_{0}]$; suppose also that we have a set
$A\subseteq \tilde{I_0}$ with the property that \beq \label{masx}
\:\exists\:\d\in(0,1)\:st\:\:\:\:\:\:\:\forall
\:P=[\a,\o,I]\in\p\:\:\operatorname{we\:have\:}\:\:|I^{*}\cap A|\leq
\d |I|. \eeq

Then  $\forall\:f\in L^2(\TT)$ we have \beq \label{cut}
\left\|\chi_{A}{T^{\p}}^{*}f\right\|_2\lesssim\d^{\frac{1}{2}}\left\|f\right\|_{2}.
\eeq
\begin{proof}

We start the proof with the following observation: if
$l_0(x)=c_0+2b_0x$ is the central line of $P_0$ then \eqref{cut} is
equivalent with
$$\left\|\chi_{A}Q_{b_0}^{*}M_{c_0}^{*}{T^{\p}}^{*}M_{c_0}Q_{b_0}f\right\|_2\lesssim\d^{\frac{1}{2}}\left\|f\right\|_{2}\:.$$

Consequently, we may assume that the top frequency line $l_{P_0}$
coincides with the real axis (indeed, for the general case, taking
as usual
${{\T}^{\p}}^{*}=Q_{b_0}^{*}M_{c_0}^{*}{T^{\p}}^{*}M_{c_0}Q_{b_0}$,
one may repeat the procedure appearing below, by using relations
\eqref{t} and \eqref{meanval} in \eqref{ets}).

Another observation is that from the structure of the intervals
$\left\{I_P^{*}\right\}$ we know that even though they are not
necessarily dyadic, each $I_P^{*}$ can be written as a union of at
most 4 dyadic intervals with the same length - call them
$\left\{I_{P,j}\right\}_{j=1}^{4}$. With this done, set
$$S=\left\{I_{P,j}\:|\:P\in\p\:\:\:\&\:\:\:j\in\left\{1,\ldots 4\right\}\right\}\:.$$

Suppose now that $\p$ is a finite collection of tiles. Also, define
$\J$ the collection of maximal dyadic intervals $I$ with the
property
$$(\ast)\:\:\:\:\:\:\:\forall\:\:J\in S\:\:\operatorname{if}\:\:J\cap I\not=\es\:\:\:\operatorname{then}\:\:\:I\subseteq
J\:.$$

Set then $\tilde{\J}$ to be any dyadic partition of $[0,1]$ that
contains $\J$. Now, by inspecting \eqref{cut}, we remark that we may
consider $A\subset supp\:{T^{\p}}^{*}$. Then from the maximality of
$\J$ and \eqref{masx}, we deduce
 \beq \label{den}
\forall\:J\in\tilde{\J}\:\:\:\:\:\:\:\:\:\:\:|A\cap J|\lesssim\d
|J|\:. \eeq

On the other hand, we also have
$$\chi_{A}{T^{\p}}^{*}f(x)=\sum_{J\in\tilde{\J}}\chi_{J\cap
A}\left\{\sum_{P\in\p}{T_{P}}^{*}f(x)\right\}\:.$$

Now our proof relies on the relation ($x\in J$ fixed, and
$J\in\tilde{\J}$) \beq \label{bmo}
\left|{T^{\p}}^{*}f(x)-\frac{1}{|J|}\int_{J}{T^{\p}}^{*}f(s)ds\right|\lesssim
\frac{1}{|J|}\int_{J}Mf(s)ds\:. \eeq

If we accept this for the moment, then, denoting
$$M_{\tilde{\J}}f(x)=\sum_{J\in\tilde{\J}}\chi_{J}(x)\sup_{J\subseteq I}\frac{1}{|I|}\int_{I}|f|(s)ds\:,$$
 we have that
$$\left|{T^{\p}}^{*}f(x)\right|\lesssim M_{\tilde{\J}}\left({T^{\p}}^{*}f\right)(x)\:+\:M_{\tilde{\J}}(Mf)(x)\:.$$

Now based on  \eqref{den}, we see that the relation \eqref{sets} is
satisfied for  $E_J=A\cap J$, and so we conclude
$$ \left|\chi_{A}{T^{\p}}^{*}f(x)\right|\lesssim M_{\d}\left({T^{\p}}^{*}f\right)(x)\:+\:M_{\d}(Mf)(x)\:,$$
which combined with \eqref{fmax} implies \eqref{cut}.

We return now at  \eqref{bmo}. For fixed $J\in\J$  and $x\in J$ we
have
$$\left|{T^{\p}}^{*}f(x)-\frac{1}{|J|}\int_{J}{T^{\p}}^{*}f(s)ds\right|=$$
\beq \label{ets}
\left|\frac{1}{|J|}\int_{J}\left\{\sum_{{P\in\p}\atop{|I_P|\geq
|J|}}\int_{\TT}\left[\f_k(x-y)-\f_k(s-y)\right]f(y)\chi_{E(P)}(y)dy\right\}ds\right|\lesssim
\eeq
$$\frac{1}{|J|}\int_{J}\left\{\sum_{2^{-k}\geq |J|}2^{k}|J|Mf(s)\right\}ds\lesssim M_{\tilde{\J}}(Mf)(x)\:,$$
and the proof of our claim is now complete.

\end{proof}
\end{l1}

\section{\bf Remarks}

1) Using interpolation methods, one can show\footnote{See also
\cite{2}, Section 8.} that the previous results can be extended to
handle the $L^p$ case ($1<p<\infty$).

2) The general polynomial phase case requires further generalization
of the tiles to curved regions in the time-frequency plane. We hope
to address this subject in the future.

\begin{comment}
It seems that the same approach works for the full conjecture
(general polynomial phase). We intend to address this subject in a
forthcoming paper.

While for the full conjecture (general polynomial phase)
 the same approach should remain effective, here the main
difference lies in quantifying the geometric interaction between
tiles. More exactly, the ``tiles" which will arise will no longer be
parallelograms associated to straight lines, but rather curved
shapes associated to higher-degree polynomial curves; in this
context, the notion of an appropriate degree of overlapping between
tiles becomes more complicated.
\end{comment}

\end{document}